\documentclass[a4paper,reqno]{amsart}

\usepackage{amssymb}

\newtheorem{theorem}[equation]{Theorem}
\newtheorem{lemma}[equation]{Lemma}
\newtheorem{corollary}[equation]{Corollary}
\newtheorem{proposition}[equation]{Proposition}

\numberwithin{equation}{section}

\theoremstyle{definition}
\newtheorem{definition}[equation]{Definition}
\newtheorem*{example*}{Example}

\newtheorem{remark}[equation]{Remark}
\newtheorem*{remark*}{Remark}

\newcommand{\bZ}{{\mathbb Z}}

\newcommand{\frg}{{\mathfrak g}}
\newcommand{\frtg}{{\tilde{\mathfrak g}}}

\newcommand{\frh}{{\mathfrak h}}
\newcommand{\frb}{{\mathfrak b}}

\newcommand{\frf}{{\mathfrak f}}
\newcommand{\fre}{{\mathfrak e}}

\newcommand{\frs}{{\mathfrak s}}

\newcommand{\frgo} {{\frg_{\bar 0}}}
\newcommand{\frguno} {{\frg_{\bar 1}}}

\newcommand{\calS}{{\mathcal S}}

\newcommand{\calJord}{{\mathcal{J}ord}}
\newcommand{\subo}{_{\bar 0}}
\newcommand{\subuno}{_{\bar 1}}

\newcommand{\bil}{{\textup{b}}}

\DeclareMathOperator{\eespan}{span}
\providecommand{\espan}[1]{\eespan\left\{ #1\right\}}

 \newcommand{\tri}{\mathfrak{tri}}

 \newcommand{\frsl}{{\mathfrak{sl}}}
 \newcommand{\frsp}{{\mathfrak{sp}}}
 \newcommand{\frso}{{\mathfrak{so}}}
 \newcommand{\frpsl}{{\mathfrak{psl}}}
 \newcommand{\frgl}{{\mathfrak{gl}}}
 \newcommand{\frpgl}{{\mathfrak{pgl}}}
 \newcommand{\frosp}{{\mathfrak{osp}}}

  \newcommand{\frbr}{{\mathfrak{br}}}
   \newcommand{\frel}{{\mathfrak{el}}}

 \DeclareMathOperator{\tr}{tr}
  \DeclareMathOperator{\str}{str}
 \DeclareMathOperator{\ad}{ad}
 \DeclareMathOperator{\der}{\mathfrak{der}}
 
 \DeclareMathOperator{\inder}{\mathfrak{inder}}
 
 \DeclareMathOperator{\End}{End}
 
 \DeclareMathOperator{\Mat}{Mat}

 \DeclareMathOperator{\charac}{char}

\newcommand{\OTS}{orthogonal triple system}

\def\bigstrut{\vrule height 12pt width 0ptdepth 2pt}

\def\hregleta{\hrule height .5pt}
\def\hreglon{\hrule height1pt}
\def\vreglon{\vrule 
    width 1pt}

\def\hreglonfill{\leaders\hreglon\hfill}

\def\hregletafill{\leaders\hregleta\hfill}

\newenvironment{romanenumerate}
 {\begin{enumerate}
 
 }{\end{enumerate}}

\begin{document}

\title{Models of some simple modular Lie superalgebras}

\author[Alberto Elduque]{Alberto Elduque$^{\star}$}
 \thanks{$^{\star}$ Supported by the Spanish Ministerio de
 Educaci\'{o}n y Ciencia
 and FEDER (MTM 2007-67884-C04-02) and by the
Diputaci\'on General de Arag\'on (Grupo de Investigaci\'on de
\'Algebra)}
 \address{Departamento de Matem\'aticas e
 Instituto Universitario de Matem\'aticas y Aplicaciones,
 Universidad de Zaragoza, 50009 Zaragoza, Spain}
 \email{elduque@unizar.es}


\date{May 9, 2008}



\begin{abstract}
Models of the exceptional simple modular Lie superalgebras in characteristic $p\geq 3$, that have appeared in the classification due to Bouarroudj, Grozman and Leites  \cite{BGLCartan} of the Lie superalgebras with indecomposable symmetrizable Cartan matrices, are provided. The models relate these exceptional Lie superalgebras to some low dimensional nonassociative algebraic systems.
\end{abstract}

\maketitle


\section*{Introduction}

The finite dimensional modular Lie superalgebras with indecomposable symmetrizable Cartan matrices over algebraically closed fields  are classified in \cite{BGLCartan} under some extra technical hypotheses. Their results assert
that, for characteristic $\geq 3$, apart from the Lie superalgebras obtained as the analogues of the Lie superalgebras in the classification in characteristic $0$ \cite{Kac-Lie}, by reducing the Cartan matrices modulo $p$, there are the following exceptions that have to be added to the list of known simple Lie superalgebras:

\begin{romanenumerate}
\item Two exceptions in characteristic $5$: $\frbr(2;5)$ and $\frel(5;5)$. (The superalgebra $\frel(5;5)$ first appeared in \cite{EldModular}.)

\item A family of exceptions given by the Lie superalgebras that appear in the \emph{Supermagic Square} in characteristic $3$ considered in \cite{CE1,CE2}. With the exception of $\frg(3,6)=\frg(S_{1,2},S_{4,2})$ these Lie superalgebras first appeared  in \cite{EldNew} and \cite{EldModular}.

\item Another two exceptions in characteristic $3$, similar to the ones in characteristic $5$: $\frbr(2;3)$ and $\frel(5;3)$.
\end{romanenumerate}

\smallskip

The Lie superalgebra $\frel(5;5)$ was shown in \cite{EldModular} to be related to Kac's $10$-dimensional exceptional Jordan superalgebra, by means of the Tits construction of Lie algebras in terms of alternative and Jordan algebras \cite{Tits66}.

The purpose of this paper is to provide models of the other three exceptions: $\frbr(2;3)$ and $\frel(5;3)$ in characteristic $3$, and $\frbr(2;5)$ in characteristic $5$.

Actually, the superalgebra $\frbr(2;3)$ already appeared in \cite[Theorem 3.2(i)]{EldNew} related to a symplectic triple system of dimension $8$. Here it will be shown to be related to a nice five dimensional orthosymplectic triple system.

The Lie superalgebra $\frel(5;3)$ will be shown to be a maximal subalgebra of the Lie superalgebra $\frg(8,3)=\frg(S_8,S_{1,2})$ in the Supermagic Square. Furthermore, it will be shown to be related to an orthogonal triple system defined on the direct sum of two copies of the octonions  and, finally, it will be proved to be the Lie superalgebra of derivations of a specific orthosymplectic triple system, and this latter result will relate $\frel(5;3)$ to the Lie superalgebra $\frg(6,6)=\frg(S_{4,2},S_{4,2})$ in the Supermagic Square.

Finally, a very explicit model of the Lie superalgebra $\frbr(2;5)$ will be constructed.

\medskip

The paper is organized as follows. The construction of the Extended Magic Square (or Supermagic Square) in characteristic $3$ in terms of composition superalgebras is recalled in \S 1. Then, in \S 2, the Lie superalgebra $\frel(5;3)$ (in characteristic $3$) is shown to be a maximal subalgebra of the Lie superalgebra $\frg(S_8,S_{1,2})$ in the Supermagic Square. This gives a very concrete realization of $\frel(5;3)$ in terms of simple components: copies of the three dimensional simple Lie algebra $\frsl_2$ and of its natural two-dimensional module. Orthogonal triple systems are reviewed in \S 3 and the Lie superalgebra $\frel(5;3)$ is shown to be isomorphic to the Lie superalgebra of an orthogonal triple system defined on the direct sum of two copies of the split Cayley algebra. Then the orthosymplectic triple systems, which extend both the orthogonal and symplectic triple systems, are recalled in \S 4. A very simple such system is defined on the set of trace zero elements of the $4\vert 2$ dimensional composition superalgebra $B(4,2)$. (The dimension being $4\vert 2$ means that the even part has dimension $4$ and the odd part dimension $2$.) The Lie superalgebra naturally attached to this orthosymplectic triple system is shown to be isomorphic to the Lie superalgebra $\frbr(2;3)$. \S 5 deals with another distinguished orthosymplectic triple system, which lives inside the Lie superalgebra $\frg(S_8,S_{1,2})$ in the Supermagic Square. It turns out that the Lie superalgebra $\frel(5;3)$ is isomorphic to the Lie superalgebra of derivations of this system. This shows also how $\frel(5;3)$ embeds in the Lie superalgebra $\frg(S_{4,2},S_{4,2})$ of the Supermagic Square. Finally, \S 6 is devoted to give an explicit model of the Lie superalgebra $\frbr(2;5)$ (in characteristic $5$) in terms of two copies of $\frsl_2$ and of their natural modules.

\medskip

All the vector spaces and superspaces considered in this paper will be assumed to be finite dimensional over a ground field $k$ of characteristic $\ne 2$. In dealing with elements of a superspace $V=V\subo\oplus V\subuno$, an expression like $(-1)^{uv}$, for homogeneous elements $u,v$, is a shorthand for $(-1)^{p(u)p(v)}$, where $p$ is the parity function.

\bigskip

\section{The Supermagic Square in characteristic $3$}

Recall that an algebra $C$ over a field $k$ is said to be a \emph{composition algebra} if it is endowed with a regular quadratic form $q$ (that is, its polar form $\bil(x,y)=q(x+y)-q(x)-q(y)$ is a nondegenerate symmetric bilinear form) such that $q(xy)=q(x)q(y)$ for any $x,y\in C$. The unital composition algebras will be termed \emph{Hurwitz algebras}. On the other hand, a composition algebra is said to be \emph{symmetric} in case the polar form is associative: $\bil(xy,z)=\bil(x,yz)$.

\smallskip

Hurwitz algebras are the well-known algebras that generalize the
classical real division algebras of the real and complex numbers,
quaternions and octonions. Over any algebraically closed field $k$,
there are exactly four of them: $k$, $k\times k$, $\Mat_2(k)$ and
$C(k)$ (the split Cayley algebra), with dimensions $1$, $2$, $4$ and
$8$.

\smallskip

Let us superize the above concepts.

A quadratic superform on a $\bZ_2$-graded vector space
$U=U\subo\oplus U\subuno$ over a field $k$ is a pair
$q=(q\subo,\bil)$ where $q\subo :U\subo\rightarrow k$ is a quadratic
form, and
 $\bil:U\times U\rightarrow k$ is a supersymmetric even bilinear form
such that $\bil\vert_{U\subo\times U\subo}$ is the polar of $q\subo$:
\[
\bil(x\subo,y\subo)=q\subo(x\subo+y\subo)-q\subo(x\subo)-q\subo(y\subo)
\]
for any $x\subo,y\subo\in U\subo$.

The quadratic superform $q=(q\subo,\bil)$ is said to be
\emph{regular} if the bilinear form $\bil$  is
nondegenerate.

\smallskip

Then a superalgebra $C=C\subo\oplus C\subuno$ over $k$, endowed with
a regular quadratic superform $q=(q\subo,\bil)$, called the
\emph{norm}, is said to be a \emph{composition superalgebra} (see
\cite{EldOkuCompoSuper}) in case
\begin{subequations}\label{eq:norm}
\begin{align}
&q\subo(x\subo y\subo)=q\subo(x\subo)q\subo(y\subo),\label{eq:qcompo1}\\
&\bil(x\subo y,x\subo z)=q\subo(x\subo)\bil(y,z)=\bil(yx\subo,zx\subo),\label{eq:qcompo2}\\
&\bil(xy,zt)+(-1)^{  x   y  +
x
 z +  y   z }\bil(zy,xt)=(-1)^{
 y   z }\bil(x,z)\bil(y,t),\label{eq:qcompo3}
\end{align}
\end{subequations}
for any $x\subo,y\subo\in C\subo$ and homogeneous elements
$x,y,z,t\in C$. Since the characteristic of the ground field is assumed to be not $2$, equation \eqref{eq:qcompo3} already implies \eqref{eq:qcompo1} and \eqref{eq:qcompo2}.

The unital composition superalgebras are termed \emph{Hurwitz
superalgebras}, while a composition superalgebra is said to be
\emph{symmetric} in case its bilinear form is associative, that is,
\[ \bil(xy,z)=\bil(x,yz),
\]
for any $x,y,z$.

\smallskip

Only over fields of characteristic $3$  there appear nontrivial
Hurwitz superalgebras (see \cite{EldOkuCompoSuper}):

\begin{itemize}

\item Let $V$ be a two dimensional vector space over a field $k$,
endowed with a nonzero alternating bilinear form $\langle .\vert
.\rangle$ (that is $\langle v\vert v\rangle =0$ for any $v\in V$).  Consider the superspace $B(1,2)$ (see \cite{She97}) with
\begin{equation}\label{eq:B12a}
B(1,2)\subo =k1,\qquad\text{and}\qquad B(1,2)\subuno= V,
\end{equation}
endowed with the supercommutative multiplication given by
\[
1x=x1=x\qquad\text{and}\qquad uv=\langle u\vert v\rangle 1
\]
for any $x\in B(1,2)$ and $u,v\in V$, and with the quadratic
superform $q=(q\subo,\bil)$ given by:
\begin{equation}\label{eq:B12b}
q\subo(1)=1,\quad \bil(u,v)=\langle u\vert v\rangle,
\end{equation}
for any $u,v\in V$. If the characteristic of $k$ is equal to $3$, then
$B(1,2)$ is a Hurwitz superalgebra (\cite[Proposition
2.7]{EldOkuCompoSuper}).

\smallskip

\item Moreover, with $V$ as before, let $f\mapsto \bar f$ be the
associated symplectic involution on $\End_k(V)$ (so $\langle
f(u)\vert v\rangle =\langle u\vert\bar f(v)\rangle$ for any $u,v\in
V$ and $f\in\End_k(V)$). Consider the superspace $B(4,2)$ (see
\cite{She97}) with
\begin{equation}\label{eq:B42}
B(4,2)\subo=\End_k(V),\qquad\text{and}\qquad B(4,2)\subuno=V,
\end{equation}
with multiplication given by the usual one (composition of maps) in
$\End_k(V)$, and by
\[
\begin{split}
&v\cdot f=f(v)=\bar f\cdot v \in V,\\
&u\cdot v=\langle .\vert u\rangle v\in \End_k(V)
\end{split}
\]
for any $f\in\End_k(V)$ and $u,v\in V$, where $\langle .\vert u\rangle v$ denotes the endomorphism $w\mapsto \langle w\vert u\rangle v$; and with quadratic superform
$q=(q\subo,\bil)$ such that
\[
q\subo(f)=\det(f),\qquad\bil(u,v)=\langle u\vert v\rangle,
\]
for any $f\in \End_k(V)$ and $u,v\in V$. If the
characteristic is equal to $3$, $B(4,2)$ is a Hurwitz superalgebra
(\cite[Proposition 2.7]{EldOkuCompoSuper}).

\end{itemize}

\smallskip

Given any Hurwitz superalgebra $C$ with norm $q=(q\subo,\bil)$, its
standard involution is given by
\[
x\mapsto \bar x=\bil(x,1)1-x.
\]
A new product can be defined on $C$ by means of
\begin{equation}\label{eq:paraHurwitz}
x\bullet y=\bar x\bar y.
\end{equation}
The resulting superalgebra, denoted by $\bar C$, is called the
\emph{para-Hurwitz superalgebra} attached to $C$, and it turns out to be a symmetric composition superalgebra.

\smallskip

Given a symmetric composition superalgebra $S$, its \emph{triality
Lie superalgebra} $\tri(S)=\tri(S)\subo\oplus\tri(S)\subuno$ is
defined by:
\begin{multline*}
\tri(S)_{\bar i}=\{ (d_0,d_1,d_2)\in\frosp(S,q)^3_{\bar i}:\\
d_0(x\bullet y)=d_1(x)\bullet y+(-1)^{i  x }x\bullet
d_2(y)\ \forall x,y\in S\subo\cup S\subuno\},
\end{multline*}
where $\bar i= \bar 0,\bar 1$, and $\frosp(S,q)$ denotes the
associated orthosymplectic Lie superalgebra. The bracket in
$\tri(S)$ is given componentwise.

Now, given two symmetric composition superalgebras $S$ and $S'$, one can form (see
\cite[\S 3]{CE1}, or \cite{ElduqueMagicI} for the non-super situation) the Lie superalgebra:
\[
\frg=\frg(S,S')=\bigl(\tri(S)\oplus\tri(S')\bigr)\oplus\bigl(\oplus_{i=0}^2
\iota_i(S\otimes S')\bigr),
\]
where $\iota_i(S\otimes S')$ is just a copy of $S\otimes S'$
($i=0,1,2$),  with bracket given by:

\begin{itemize}
\item the Lie bracket in $\tri(S)\oplus\tri(S')$, which thus becomes  a Lie subalgebra of $\frg$,
\smallskip

\item $[(d_0,d_1,d_2),\iota_i(x\otimes
 x')]=\iota_i\bigl(d_i(x)\otimes x'\bigr)$,
\smallskip

\item
 $[(d_0',d_1',d_2'),\iota_i(x\otimes
 x')]=(-1)^{  d_i'   x }\iota_i\bigl(x\otimes d_i'(x')\bigr)$,
\smallskip

\item $[\iota_i(x\otimes x'),\iota_{i+1}(y\otimes y')]=(-1)^{
x'   y }
 \iota_{i+2}\bigl((x\bullet y)\otimes (x'\bullet y')\bigr)$ (indices modulo
 $3$),
\smallskip

\item $[\iota_i(x\otimes x'),\iota_i(y\otimes y')]=
 (-1)^{  x   x' +  x   y'  +
   y   y' }
 \bil'(x',y')\theta^i(t_{x,y})$ \newline \null\hspace{2.5 in} $+
 (-1)^{  y   x' }
 \bil(x,y)\theta'^i(t'_{x',y'})$,

\end{itemize}
for any $i=0,1,2$ and homogeneous $x,y\in S$, $x',y'\in S'$,
$(d_0,d_1,d_2)\in\tri(S)$, and $(d_0',d_1',d_2')\in\tri(S')$. Here
$\theta$ denotes the natural automorphism
$\theta:(d_0,d_1,d_2)\mapsto (d_2,d_0,d_1)$ in $\tri(S)$, while $t_{x,y}$ is defined by
\begin{equation}\label{eq:txy}
t_{x,y}=\bigl(\sigma_{x,y},\tfrac{1}{2}\bil(x,y)1-r_xl_y,\tfrac{1}{2}\bil(x,y)1-l_xr_y\bigr)
\end{equation}
with $l_x(y)=x\bullet y$, $r_x(y)=(-1)^{xy}y\bullet x$, and
\begin{equation}\label{eq:sigmaxy}
\sigma_{x,y}(z)=(-1)^{yz}\bil(x,z)y-(-1)^{x(y+z)}\bil(y,z)x
\end{equation}
for homogeneous $x,y,z\in S$. Also $\theta'$
and $t'_{x',y'}$ denote the analogous elements for $\tri(S')$.

\smallskip

Over a field $k$ of characteristic $3$, let $S_r$ ($r=1$, $2$, $4$ or $8$) denote the para-Hurwitz algebra attached to the split Hurwitz
algebra of dimension $r$ (this latter algebra being either $k$, $k\times k$,
$\Mat_2(k)$ or $C(k)$). Also, denote by
$S_{1,2}$ the para-Hurwitz superalgebra
$\overline{B(1,2)}$, and by $S_{4,2}$ the
para-Hurwitz superalgebra $\overline{B(4,2)}$. Then the Lie superalgebras $\frg(S,S')$, where $S,S'$ run over $\{S_1,S_2,S_4,S_8,S_{1,2},S_{4,2}\}$, appear in Table \ref{ta:supermagicsquare}, which has been obtained in \cite{CE1}.

\begin{table}[h!]
$$
\vbox{\offinterlineskip
 \halign{\hfil\ $#$\ \hfil&%
 \vreglon #%
 &\hfil\ $#$\ \hfil&\hfil\ $#$\ \hfil
 &\hfil\ $#$\ \hfil&\hfil\ $#$\ \hfil&%
 \vrule  depth 4pt width .5pt #%
 &\hfil\ $#$\ \hfil&\hfil\ $#$\ \hfil\cr
 \bigstrut &&S_1&S_2&S_4&S_8&\omit%
    \vrule height 8pt depth 4pt width .5pt&S_{1,2}&S_{4,2}\cr
 \multispan9{\hreglonfill}\cr
 S_1&&\frsl_2&\frpgl_3&\frsp_6&\frf_4&&\frpsl_{2,2}&\frsp_6\oplus (14)\cr
 \bigstrut S_2&& &\omit$\frpgl_3\oplus \frpgl_3$&\frpgl_6&\tilde \fre_6&%
     &\bigl(\frpgl_3\oplus\frsl_2\bigr)\oplus\bigl(\frpsl_3\otimes (2)\bigr)&
    \frpgl_6\oplus (20)\cr
 \bigstrut S_4&& & &\frso_{12}&\fre_7&
   &\bigl(\frsp_6\oplus\frsl_2\bigr)\oplus\bigl((13)\otimes (2)\bigr)
    &\frso_{12}\oplus spin_{12}\cr
 \bigstrut S_8&& & & &\fre_8&
    &\bigl(\frf_4\oplus\frsl_2\bigr)\oplus\bigl((25)\otimes (2)\bigr)&
      \fre_7\oplus (56)\cr
 \multispan9{\hregletafill}\cr
 \bigstrut S_{1,2}&& & & & & & \frso_7\oplus 2spin_7 &\frsp_8\oplus(40)\cr
 \bigstrut S_{4,2}&& & & & & & & \frso_{13}\oplus spin_{13}\cr}}
$$
\bigskip
\caption{Supermagic Square (characteristic
$3$)}\label{ta:supermagicsquare}
\end{table}

Since the construction of $\frg(S,S')$ is symmetric, only the
entries above the diagonal are needed. In Table
\ref{ta:supermagicsquare}, $\frf_4,\fre_6,\fre_7,\fre_8$ denote the
simple exceptional classical Lie algebras, $\tilde\fre_6$ denotes a
$78$ dimensional Lie algebras whose derived Lie algebra is the $77$
dimensional simple Lie algebra $\fre_6$ in characteristic
$3$. The even and odd parts of the nontrivial superalgebras in the
table which have no counterpart in the classification in
characteristic $0$ (\cite{Kac-Lie}) are displayed, $spin$ denotes the
spin module for the corresponding orthogonal Lie algebra, while
$(n)$ denotes a module of dimension $n$, whose precise description is given in \cite{CE1}. Thus, for example,
$\frg(S_4,S_{1,2})$ is a Lie superalgebra whose even part is
(isomorphic to) the direct sum of the symplectic Lie algebra
$\frsp_6$ and of $\frsl_2$, while its odd part is the tensor
product of a $13$ dimensional module for $\frsp_6$ and the
natural $2$ dimensional module for $\frsl_2$.

In Table \ref{ta:supermagicsquare2}, a more precise description of the Lie superalgebras that appear in the Supermagic Square is given. This table displays the even parts and the highest weights of the odd parts.
The numbering of the roots follows Bourbaki's conventions \cite{Bourbaki}. The fundamental dominant weight for $\frsl_2$ will be denoted by $\omega$, while the fundamental dominant weights for a Lie algebra with a Cartan matrix of order $n$ will be denoted by $\omega_1,\ldots,\omega_n$. Below each entry, there appears the result in \cite{CE1} where the result can be found.

\def\mystrut{\vrule height 20pt width 0pt depth 8pt}
\begin{table}[h!]
$$
\vbox{\offinterlineskip
 \halign{\hfil\ $#$\ \hfil&%
 \vreglon #%
 &\hfil\ $#$\ \hfil&\hfil\ $#$\ \hfil\cr
  &depth 4pt&
    S_{1,2}&S_{4,2}\cr
 \multispan4{\hreglonfill}\cr
 \mystrut S_1&&
    \frpsl_{2,2}&\frsp_6\oplus V(\omega_3)\cr
    &&& \textrm{{\tiny\cite[Proposition 5.3]{CE1}}}\cr
 \mystrut S_2&&
    \bigl(\frpgl_3\oplus\frsl_2\bigr)\oplus\bigl(V(\omega_1+\omega_2)\otimes V(\omega)\bigr)&
    \frpgl_6\oplus V(\omega_3))\cr
    &&\textrm{{\tiny\cite[Proposition 5.28]{CE1}}}
        &\textrm{{\tiny\cite[Proposition 5.16]{CE1}}}\cr
 \mystrut S_4&&
      \bigl(\frsp_6\oplus\frsl_2\bigr)\oplus\bigl(V(\omega_2)\otimes V(\omega)\bigr)
    &\frso_{12}\oplus V(\omega_6)\cr
    &&\textrm{{\tiny\cite[Proposition 5.24]{CE1}}}
    &\textrm{{\tiny\cite[Proposition 5.5]{CE1}}}\cr
 \mystrut S_8&&
    \bigl(\frf_4\oplus\frsl_2\bigr)\oplus\bigl(V(\omega_4)\otimes V(\omega)\bigr)&
      \fre_7\oplus V(\omega_7)\cr
    \vrule width0pt depth 10pt&&\textrm{{\tiny\cite[Proposition 5.26]{CE1}}}
        &\textrm{{\tiny\cite[Proposition 5.8]{CE1}}}\cr
 \multispan4{\hregletafill}\cr
 \mystrut S_{1,2}&&
       \frso_7\oplus 2V(\omega_3) &\frsp_8\oplus V(\omega_3)\cr
       &&\textrm{{\tiny\cite[Proposition 5.19]{CE1}}}
          &\textrm{{\tiny\cite[Proposition 5.12]{CE1}}}\cr
 \mystrut S_{4,2}&&
         & \frso_{13}\oplus V(\omega_6)\cr
         &&&\textrm{{\tiny\cite[Proposition 5.10]{CE1}}}\cr}}
$$
\bigskip
\caption{Even and odd parts in the Supermagic Square}\label{ta:supermagicsquare2}
\end{table}

\bigskip

A precise description of these modules and of the superalgebras in Table \ref{ta:supermagicsquare2}
as Lie superalgebras with a Cartan matrix is given in \cite{CE1}. All the
inequivalent Cartan matrices for these simple Lie superalgebras are listed in \cite{BGL}.

With the exception of $\frg(S_{1,2},S_{4,2})$, all these superalgebras have appeared previously in \cite{EldNew} and \cite{EldModular}. Some relationships between the Lie superalgebras $\frg(S_{1,2},S)$ and $\frg(S_{4,2},S)$ to other algebraic structures have been considered in \cite{CE2}.

\bigskip

\section{The Lie superalgebra $\frel(5;3)$}

The aim of this section is to show how the Lie superalgebra $\frel(5;3)$ embeds in a nice way as a maximal subalgebra in the simple Lie superalgebra $\frg(S_8,S_{1,2})$ of the Supermagic Square.

Throughout this section the characteristic of the ground field $k$ will be assumed to be $3$.

The para-Hurwitz superalgebra $S_{1,2}=\overline{B(1,2)}$ is described as $S_{1,2}=k1\oplus V$ (see \eqref{eq:B12a} and \eqref{eq:B12b}), where $(S_{1,2})\subo=k1$ is a copy of the ground field, and $(S_{1,2})\subuno=V$ is a two dimensional vector space equipped with a nonzero alternating bilinear form $\langle.\vert.\rangle$. The multiplication is given by:
\begin{equation}\label{eq:S12product}
1\bullet 1=1,\quad 1\bullet u=-u=u\bullet 1,\quad u\bullet v=\langle u\vert v\rangle 1,
\end{equation}
for any $u,v\in V$, and the norm $q=(q\subo,\bil)$ is given by
\begin{equation}\label{eq:S12bil}
q\subo(1)=1,\qquad \bil(u,v)=\langle u\vert v\rangle,
\end{equation}
for any $u,v\in V$.

Recall from \cite{EldOkuCompoSuper} or \cite[Corollary 2.12]{CE1} that the triality Lie superalgebra of $S_{1,2}$ is given by:
\begin{equation}\label{eq:triS12}
\tri(S_{1,2})=\{(d,d,d): d\in \frosp(S_{1,2},q)\},
\end{equation}
and thus $\tri(S_{1,2})$ can (and will) be identified with the Lie superalgebra $\frb_{0,1}=\frsp(V)\oplus V$ (see \cite[(2.18)]{CE1}), with even part $\frsp(V)\,(\cong \frsl_2)$, odd part $V$, where $[\rho,v]=\rho(v)$ and $[u,v]=\gamma_{u,v}$ for any $\rho\in \frsp(V)$ and $u,v\in V$, with $\gamma_{u,v}=\langle u\vert .\rangle v+\langle v\vert .\rangle u$.

Besides, the action of $\frb_{0,1}$ on $S_{1,2}$ is given by:
\[
\rho:\begin{cases} 1\mapsto 0&\\ u\mapsto \rho(u)&\end{cases}\qquad\qquad
u:\begin{cases} 1\mapsto -u&\\ v\mapsto -\langle u\vert v\rangle 1,&\end{cases}
\]
for any $\rho\in \frsp(V)$ and $u,v\in V$ (see \cite[(2.16)]{CE1}).

\smallskip

Consider now the Lie superalgebra $\frg(S_8,S_{1,2})$ in the Supermagic Square:
\[
\frg(S_8,S_{1,2})=\bigl(\tri(S_8)\oplus\tri(S_{1,2})\bigr)\oplus
\bigl(\oplus_{i=0}^2\iota_i(S_8\otimes S_{1,2})\bigr).
\]
This is $\bZ_2\times\bZ_2$-graded with
\[
\begin{split}
\frg(S_8,S_{1,2})_{(0,0)}&=\tri(S_8)\oplus\tri(S_{1,2}),\\
\frg(S_8,S_{1,2})_{(1,0)}&=\iota_0(S_8\otimes S_{1,2}),\\
\frg(S_8,S_{1,2})_{(0,1)}&=\iota_1(S_8\otimes S_{1,2}),\\
\frg(S_8,S_{1,2})_{(1,1)}&=\iota_2(S_8\otimes S_{1,2}),
\end{split}
\]
and, therefore, the linear map $\tau$, defined by
\[
\tau =\begin{cases} id&\text{on $\frg(S_8,S_{1,2})_{(0,0)}\oplus \frg(S_8,S_{1,2})_{(1,0)}$,}\\
 -id&\text{on $\frg(S_8,S_{1,2})_{(0,1)}\oplus\frg(S_8,S_{1,2})_{(1,1)}$,}
\end{cases}
\]
is a Lie superalgebra automorphism. On the other hand, the grading automorphism $\sigma$:
\[
\sigma=\begin{cases} id&\text{on $\frg(S_8,S_{1,2})\subo$,}\\
 -id&\text{on $\frg(S_8,S_{1,2})\subuno$,}
\end{cases}
\]
commutes with $\tau$. Consider the order two automorphism $\xi=\sigma\tau=\tau\sigma$, which provides a $\bZ_2$-grading of $\frg(S_8,S_{1,2})$ with even and odd components given by:
\begin{equation}\label{eq:g8301}
\begin{split}
\frg(S_8,S_{1,2})_+&=\bigl(\tri(S_8)\oplus\frsp(V)\bigr)\oplus \iota_0(S_8\otimes 1)\oplus\iota_1(S_8\otimes V)\oplus \iota_2(S_8\otimes V),\\
\frg(S_8,S_{1,2})_-&=V\oplus \iota_0(S_8\otimes V)\oplus\iota_1(S_8\otimes 1)\oplus\iota_2(S_8\otimes 1).
\end{split}
\end{equation}

\begin{theorem}\label{th:el53maximalg83}
In the situation above, the  subalgebra $\frg(S_8,S_{1,2})_+$ of $\frg(S_8,S_{1,2})$ fixed by the automorphism $\xi$ is a maximal subalgebra of $\frg(S_8,S_{1,2})$ isomorphic to the Lie superalgebra $\frel(5;3)$.
\end{theorem}
\begin{proof}
As a module for the subalgebra $\tri(S_8)\oplus \frsp(V)$ of $\frg(S_8,S_{1,2})_+$, the odd component $\frg(S_8,S_{1,2})_-$ relative to the $\bZ_2$-grading given by $\xi$ decomposes as the direct sum of the nonisomorphic irreducible modules:
\[
V,\quad \iota_0(S_8\otimes V),\quad \iota_1(S_8\otimes 1),\quad \iota_2(S_8\otimes 1).
\]
Actually, identifying $\tri(S_8)$ to the orthogonal Lie algebra $\frso_8$ through the projection onto the first component (this is possible because of the Local Principle of Triality \cite[\S 35]{KMRT}), $\iota_1(S_8\otimes 1)$ and $\iota_2(S_8\otimes 1)$ are the two half-spin representations of $\frso_8$, while $\iota_0(S_8\otimes V)$ is the tensor product of the natural modules for $\frso_8$ and for $\frsp(V)$, so these four modules are indeed nonisomorphic. Therefore, any $\frg(S_8,S_{1,2})_+$-submodule of $\frg(S_8,S_{1,2})_-$ is a direct sum of some of them. But the definition of the Lie bracket in $\frg(S_8,S_{1,2})$ shows that any of these spaces generates $\frg(S_8,S_{1,2})_-$ as a module over $\frg(S_8,S_{1,2})_+$. Hence $\frg(S_8,S_{1,2})_-$ is an irreducible module for $\frg(S_8,S_{1,2})_+$ and, therefore, $\frg(S_8,S_{1,2})_+$ is a maximal subalgebra of $\frg(S_8,S_{1,2})$.

From now on, the proof relies heavily on the description of $\frg(S_8,S_{1,2})$ given in
\cite[\S 5.10]{CE1} (which follows the ideas in \cite{ElduqueMagicII}). This description is obtained in terms of five vector spaces of dimension $2$: $V_1,\ldots,V_5$, endowed with nonzero alternating bilinear forms:
\begin{equation}\label{eq:g83sigmas}
\frg(S_8,S_{1,2})=\oplus_{\sigma\in \calS_{8,3}}V(\sigma),
\end{equation}
with
\[
\begin{split}
\calS_{8,3}=\bigl\{\emptyset,&\{1,2,3,4\},
        \{5\},\{1,2\},\{3,4\},\{1,2,5\},\{3,4,5\},\\
       &\{2,3\},\{1,4\},\{2,3,5\},\{1,4,5\},
       \{1,3\},\{2,4\},\{1,3,5\},\{2,4,5\}\bigr\}.
\end{split}
\]
Here $V(\emptyset)=\oplus_{i=1}^5\frsp(V_i)$, while for $\emptyset\ne\sigma=\{i_1,\ldots,i_r\}$, $V(\sigma)=V_{i_1}\otimes\cdots\otimes V_{i_r}$.  Also, any $\sigma\subseteq \{1,2,3,4,5\}$ can be thought of as an element in $\bZ_2^5$ (for instance, $\{1,3,5\}=(\bar 1,\bar 0,\bar 1,\bar 0,\bar 1)\in\bZ_2^5$), so it makes sense to consider $\sigma+\tau$ for $\sigma,\tau\subseteq \{1,2,3,4,5\}$.

The brackets $V(\sigma)\times V(\tau)\rightarrow V(\sigma+\tau)$ are nonzero scalar multiples of the `contraction maps' $\varphi_{\sigma,\tau}$ in \cite[(4.9)]{CE1}.
Under this description,
\begin{equation}\label{eq:g830sigmas}
\frg(S_8,S_{1,2})_+=\oplus_{\sigma\in\tilde\calS_{8,3}}V(\sigma),
\end{equation}
with
\[
\tilde\calS_{8,3}=\bigl\{\emptyset,\{1,2,3,4\},
       \{1,2\},\{3,4\},
       \{2,3,5\},\{1,4,5\},\{1,3,5\},\{2,4,5\}\bigr\}.
\]
Thus, the even and odd degrees (same notation as in \cite[\S 5]{CE1}) are:
\[
\begin{split}
\Phi\subo&=\{\pm 2\epsilon_i: 1\leq i\leq 5\}\cup
  \{\pm\epsilon_1\pm\epsilon_2\pm\epsilon_3\pm\epsilon_4\}\cup
    \{\pm\epsilon_i\pm\epsilon_j:(i,j)\in\{(1,2),(3,4)\}\},\\
\Phi\subuno&=\{\pm\epsilon_5\}\cup\{\pm\epsilon_i\pm\epsilon_j\pm\epsilon_5:
          (i,j)\in\{(2,3),(1,4),(1,3),(2,4)\}\}.
\end{split}
\]
With the lexicographic order given by
$0<\epsilon_1<\epsilon_2<\epsilon_3<\epsilon_4<\epsilon_5$ in \cite[\S 5.10]{CE1}, the set of irreducible degrees is
\[
\Pi=\{\alpha_1=\epsilon_5-\epsilon_2-\epsilon_4,
  \alpha_2=\epsilon_2-\epsilon_1,\alpha_3=2\epsilon_1,
 \alpha_4=\epsilon_4-\epsilon_1-\epsilon_2-\epsilon_3,\alpha_5=2\epsilon_3\},
\]
which is a $\bZ$-linearly independent set with $\Phi=\Phi\subo\cup\Phi\subuno\subseteq
\bZ\Pi$. The associated Cartan matrix is:
\[
\begin{pmatrix}
 0&-2&0&0&0\\ -1&2&-2&0&0\\ 0&-1&2&-1&0\\ 0&0&-1&2&-1\\ 0&0&0&-1&2
\end{pmatrix}
\]
which is equal to the second matrix in \cite[\S13.2]{BGLCartan} for $\frel(5;3)$ with the third and fourth rows and columns permuted. This shows that $\frg(S_8,S_{1,2})_+$ is isomorphic to $\frel(5;3)$. Note that the $4\times 4$ submatrix on the lower right corner is the Cartan matrix of type $B_4$, and indeed it corresponds to the subalgebra $\tri(S_8)\oplus\iota_0(S_8\otimes 1)$, which is isomorphic to the orthogonal Lie algebra $\frso_9$ ($\tri(S_8)$ being isomorphic to $\frso_8$ and $S_8$ to its natural module).
\end{proof}

\bigskip

\section{Orthogonal triple systems and the Lie superalgebra $\frel(5;3)$}

This section is devoted to the proof of the fact that the Lie superalgebra $\frel(5;3)$ is the Lie superalgebra associated to a particular orthogonal triple system defined on the direct sum of two copies of the split octonions.

\smallskip

Orthogonal triple systems were introduced in \cite{OkuOTS}:

\begin{definition}\label{df:OTS}
Let $T$ be a vector space over a field $k$ endowed with a nonzero
symmetric bilinear form $(.\vert.):T\times T\rightarrow k$, and a
triple product $T\times T\times T\rightarrow T$: $(x,y,z)\mapsto
[xyz]$. Then $\bigl(T,[...],(.\vert.)\bigr)$ is said to be an \emph{orthogonal triple system} if
it satisfies the following identitities:
\begin{subequations}\label{eq:OTS}
\begin{align}
&[xxy]=0\label{eq:OTSa}\\
&[xyy]=(x\vert y)y-(y\vert y)x\label{eq:OTSb}\\
&[xy[uvw]]=[[xyu]vw]+[u[xyv]w]+[uv[xyw]]\label{eq:OTSc}\\
&([xyu]\vert v)+(u\vert [xyv])=0\label{eq:OTSd}
\end{align}
\end{subequations}
for any elements $x,y,u,v,w\in T$.
\end{definition}

Equation \eqref{eq:OTSc} shows that $\inder T=\espan{[xy.]: x,y\in
T}$ is a subalgebra (actually an ideal) of the Lie algebra $\der T$
of derivations of $T$. The elements in $\inder T$ are called \emph{inner
derivations}. Because of \eqref{eq:OTSb}, if $\dim T\geq 2$, then
$\der T$ is contained in the orthogonal Lie algebra
$\frso\bigl(T,(.\vert .)\bigr)$. Also note that \eqref{eq:OTSd} is a consequence of \eqref{eq:OTSb} and \eqref{eq:OTSc} (see the comments in \cite{EldNew} after Definition 4.1).

An ideal of an orthogonal triple system $\bigl(T,[...],(.\vert .)\bigr)$ is a subspace $I$ such that $[ITT]+[TIT]+[TTI]$ is contained in $I$. The orthogonal triple system is said to be simple if it does not contain any proper ideal.

\smallskip

Some of the main properties of these systems are summarized in the next result, taken from \cite[Proposition 4.4, Theorem 4.5 and Theorem 5.1]{EldNew} (see also \cite[Theorem 4.3]{CE2}):

\begin{proposition}\label{pr:propertiesOTS}
Let $\bigl(T,[...],(.\vert.)\bigr)$ be an orthogonal triple system of dimension $\geq 2$. Then:
\begin{enumerate}
\item $\bigl(T,[...],(.\vert.)\bigr)$ is simple if and only if $(.\vert.)$ is nondegenerate.

\smallskip

\item Let $(V,\langle.\vert.\rangle)$ be a two dimensional vector space endowed with a nonzero alternating bilinear form. Let $\frs$
be a Lie subalgebra of $\der T$ containing $\inder T$. Define the
superalgebra $\frg=\frg(T,\frs)=\frgo\oplus\frguno$ with
\[
\begin{cases}
\frgo= \frsp(V)\oplus \frs&\\
\frguno=V\otimes T\,,
\end{cases}
\]
and superanticommutative multiplication given by:
\begin{itemize}
\item
the multiplication on $\frgo$ coincides with its bracket as a Lie algebra (the direct sum of the ideals $\frsp(V)$ and $\frs$);
\item
$\frgo$ acts naturally on $\frguno$, that is,
\[
[s,v\otimes x]=s(v)\otimes x,\qquad [d,v\otimes x]=v\otimes d(x),
\]
for any $s\in \frsp(V)$, $d\in\frs$, $v\in V$, and $x\in T$;
\item
for any $u,v\in V$ and $x,y\in T$:
\begin{equation}\label{eq:gToproduct}
[u\otimes x,v\otimes y]=
  -(x\vert y)\gamma_{u,v} +\langle u\vert v\rangle d_{x,y}
\end{equation}
where $\gamma_{u,v}=\langle u\vert .\rangle v+\langle v\vert
.\rangle u$ and $d_{x,y}=[xy.]$.
\end{itemize}
Then $\frg(T,\frs)$ is a Lie superalgebra. Moreover, $\frg(T,\frs)$
is simple if and only if $\frs$ coincides with $\inder T$ and $T$ is
a simple orthogonal triple system.

\smallskip

Conversely, given a  Lie superalgebra $\frg=\frgo\oplus\frguno$ with
\[
\begin{cases}
\frgo=\frsp(V)\oplus \frs &\text{(direct sum of ideals),}\\
\frguno=V\otimes T&\text{(as a module for $\frgo$),}
\end{cases}
\]
where $T$ is a module over $\frs$, by $\frsp(V)$-invariance of the
Lie bracket, equation \eqref{eq:gToproduct} is satisfied for a
symmetric bilinear form $(.\vert .):T\times T\rightarrow k$ and an
antisymmetric bilinear map $d_{.,.}:T\times T\rightarrow \frs$.
Then, if $(.\vert .)$ is not $0$ and a triple product on $T$ is
defined by means of $[xyz]=d_{x,y}(z)$, $T$ becomes and orthogonal
triple system and the image of $\frs$ in $\frgl(T)$ under the given
representation is a subalgebra of $\der T$ containing $\inder T$.

\smallskip

\item If the characteristic of the ground field $k$ is equal to $3$, define the $\bZ_2$-graded algebra
$\frtg=\frtg(T)=\frtg\subo\oplus\frtg\subuno$, with:
\[
\frtg\subo=\inder(T),\qquad\frtg\subuno=T,
\]
and anticommutative multiplication given by:
\begin{itemize}
\item
the multiplication on $\frtg\subo$ coincides with its bracket as a Lie algebra;
\item
$\frtg\subo$ acts naturally on $\frtg\subuno$, that is,
$[d,x]=d(x)$ for any $d\in\inder(T)$ and $x\in T$;
\item
$[x,y]=d_{x,y}=[xy.]$, for any $x,y\in T$.
\end{itemize}
Then $\frtg(T)$ is a Lie algebra. Moreover, $T$ is a simple orthogonal triple system if and
only if  $\frtg(T)$ is a simple $\bZ_2$-graded Lie algebra.

\end{enumerate}
\end{proposition}

The Lie superalgebra $\frg(T)=\frg(T,\inder(T))$ in item 2) above will be called the \emph{Lie superalgebra of the orthogonal triple system $T$} and, if the characteristic is $3$, the Lie algebra $\frtg(T)$ will be called the \emph{Lie algebra of the orthogonal triple system $T$}.

\smallskip

The classification of the simple finite dimensional orthogonal triple systems in characteristic $0$ appears in \cite[Theorem 4.7]{EldNew}. In characteristic $3$, there appears at least one new family of simple orthogonal triple systems, which are attached to degree $3$ Jordan algebras (see \cite[Examples 4.20]{EldNew}):

Let $J=\calJord(n,1)$ be the Jordan algebra of a nondegenerate
cubic form $n$ with basepoint $1$, over a field $k$ of
characteristic $3$, and assume that $\dim_kJ\geq 3$. Then any
$x\in J$ satisfies a cubic equation \cite[II.4]{McCrimmon}
\begin{equation}\label{eq:degree3}
x^{\circ 3}-t(x)x^{\circ 2}+s(x)x-n(x)1=0,
\end{equation}
where $t$ is its trace linear form, $s(x)$ is the spur
quadratic form and the multiplication in $J$ is denoted by $x\circ
y$. For our purposes it is enough to consider the Jordan algebras in \eqref{eq:H3Z2Z2} below.

Let $J_0=\{x\in J:t(x)=0\}$ be the subspace of trace zero
elements. Since $\charac k=3$, $t(1)=0$, so that $k\, 1\in J_0$.
Consider the quotient space $\hat J=J_0/k\, 1$. For any $x\in
J_0$, we have $s(x)=-\tfrac{1}{2}t(x^{\circ 2})$ and, by linearization of
\eqref{eq:degree3}, we get that for any $x,y\in J_0$:
\begin{equation}\label{eq:yyx}
\begin{split}
y^{\circ 2}\circ x-(x\circ y)\circ y&\equiv
 -2t(x,y)y-t(y,y)x \mod k\,1,\\
 &\equiv t(x,y)y-t(y,y)x \mod k\, 1.
\end{split}
\end{equation}
Let us denote by $\hat x$ the class of $x$ modulo $k\,1$. Since
$J_0$ is the orthogonal complement of $k\,1$ relative to the trace bilinear
form $t(a,b)=t(a\circ b)$, $t$ induces a nondegenerate symmetric
bilinear form on $\hat J$ defined by $t(\hat x,\hat y)=t(x,y)$ for
any $x,y\in J_0$. Now, for any $x,y\in J_0$ consider the inner
derivation of $J$ given by $D_{x,y}:z\mapsto x\circ (y\circ
z)-y\circ(x\circ z)$ (see \cite{JacobsonJordan}). Since the trace
form is invariant under the Lie algebra of derivations, $D_{x,y}$
leaves $J_0$ invariant, and obviously satisfies $D_{x,y}(1)=0$, so
it induces a map $d_{x,y}:\hat J\rightarrow \hat J$, $\hat
z\mapsto \widehat{D_{x,y}(z)}$ and a well defined bilinear map
$\hat J\times \hat J\rightarrow \frgl(\hat J)$, $(\hat x,\hat
y)\mapsto d_{x,y}$. Consider now the triple  product $[...]$ on
$\hat J$ defined by
\[
[\hat x\hat y\hat z]=d_{x,y}(\hat z)
\]
for any $x,y,z\in J_0$. This is well defined and satisfies
\eqref{eq:OTSa}, because of the antisymmetry of $d_{.,.}$. Also,
\eqref{eq:yyx} implies that
\[
\begin{split}
[\hat x\hat y\hat y]=d_{x,y}(\hat y)&=t(x,y)\hat y-t(y,y)\hat x\\
&=t(\hat x,\hat y)\hat y-t(\hat y,\hat y)\hat x,
\end{split}
\]
so \eqref{eq:OTSb} is satisfied too, relative to the trace
bilinear form. Since $D_{x,y}$ is a derivation of $J$ for any
$x,y\in J$, \eqref{eq:OTSc} follows immediately, while
\eqref{eq:OTSd} is a consequence of $D_{x,y}$ being a derivation and the trace $t$ being associative.

Therefore, by nondegeneracy of the trace form,
$\bigl(\hat J,[...], t(.,.)\bigr)$ is a simple \OTS\ over
$k$ \cite[Examples 4.20]{EldNew}.

Now, let $e\ne 0,1$ be an idempotent ($e^{\circ 2}=e$) of such a Jordan algebra. Changing $e$ by $1-e$ if necessary, it can be assumed that $t(e)=1$. Consider the Peirce $1$-space:
\[
J_1(e)=\{x\in J: e\circ x=\frac{1}{2}x\}
\]
Note that $J_1(e)$ is contained in $J_0$, because for any $x\in J_1(e)$, we have
\[
t(x)=2t(e\circ x)=2t((e\circ e)\circ x))=2t(e\circ(e\circ x))=\frac{1}{2}t(x),
\]
so $t(x)=0$, and since $1\in J_0(e)\oplus J_2(e)$, $J_1(e)$ embeds in $\hat J=J_0/k1$. Besides, since $J_1(e)\circ J_1(e)\subseteq J_0(e)\oplus J_2(e)$, and $(J_0(e)\oplus J_2(e))\circ J_1(e)\subseteq J_1(e)$ (see \cite[II.8]{McCrimmon}), it follows that $J_1(e)$ is an orthogonal triple subsystem of the orthogonal triple system $\hat J$ above.

In particular, let $C$ be a Hurwitz algebra over the field $k$ of characteristic $3$ with norm $q$ and polar form $\bil$, and consider the Jordan algebra $J=H_3(C,*)$ of hermitian $3\times 3$ matrices (where $(a_{ij})^*=(\bar a_{ji})$) under the symmetrized product $x\circ y=\frac{1}{2}(xy+yx)$. Let $S$ be the associated para-Hurwitz algebra. Then,
\begin{equation}\label{eq:H3Z2Z2}
\begin{split}
J=H_3(C,*)&=\left\{ \begin{pmatrix} \alpha_0 &\bar a_2& a_1\\
  a_2&\alpha_1&\bar a_0\\ \bar a_1&a_0&\alpha_2\end{pmatrix} :
  \alpha_0,\alpha_1,\alpha_2\in k,\ a_0,a_1,a_2\in S\right\}\\[6pt]
 &=\bigl(\oplus_{i=0}^2 ke_i\bigr)\oplus
     \bigl(\oplus_{i=0}^2\iota_i(S)\bigr),
\end{split}
\end{equation}
where
\begin{equation}\label{eq:eisiotas}
\begin{aligned}
e_0&= \begin{pmatrix} 1&0&0\\ 0&0&0\\ 0&0&0\end{pmatrix}, &
 e_1&=\begin{pmatrix} 0&0&0\\ 0&1&0\\ 0&0&0\end{pmatrix}, &
 e_2&= \begin{pmatrix} 0&0&0\\ 0&0&0\\ 0&0&1\end{pmatrix}, \\
 \iota_0(a)&=2\begin{pmatrix} 0&0&0\\ 0&0&\bar a\\
 0&a&0\end{pmatrix},&
 \iota_1(a)&=2\begin{pmatrix} 0&0&a\\ 0&0&0\\
 \bar a&0&0\end{pmatrix},&
 \iota_2(a)&=2\begin{pmatrix} 0&\bar a&0\\ a&0&0\\
 0&0&0\end{pmatrix},
\end{aligned}
\end{equation}
for any $a\in S$. Then $J$ is the Jordan algebra of the nondegenerate cubic form $n$ with basepoint $1$, where
\[
n(x)=\alpha_0\alpha_1\alpha_2+\bil(a_0a_1a_2,1)-\sum_{i=0}^2\alpha_iq(a_i),
\]
for $x=\sum_{i=0}^2\alpha_ie_i + \sum_{i=0}^2\iota_i(a_i)$. Here the trace form $t$ is the usual trace: $t(x)=\sum_{i=0}^2\alpha_i$.

Identify $ke_0\oplus ke_1\oplus ke_2$ with $k^3$ by
means of $\alpha_0e_0+\alpha_1e_1+\alpha_2e_2\simeq
(\alpha_0,\alpha_1,\alpha_2)$. Then the Jordan product becomes:
\begin{equation}\label{eq:Jniceproduct}
\left\{\begin{aligned}
 &(\alpha_0,\alpha_1,\alpha_2)\circ(\beta_1,\beta_2,\beta_3)=
    (\alpha_0\beta_0,\alpha_1\beta_1,\alpha_2\beta_2),\\
 &(\alpha_0,\alpha_1,\alpha_2)\circ \iota_i(a)
  =\frac{1}{2}(\alpha_{i+1}+\alpha_{i+2})\iota_i(a),\\
 &\iota_i(a)\circ\iota_{i+1}(b)=\iota_{i+2}(a\bullet b),\\
 &\iota_i(a)\circ\iota_i(b)=2\bil(a,b)\bigl(e_{i+1}+e_{i+2}\bigr),
\end{aligned}\right.
\end{equation}
for any $\alpha_i,\beta_i\in k$, $a,b\in S$, where $i=0,1,2$, and
indices are taken modulo $3$.

Now, $e=e_0$ is an idempotent of trace $1$ and the Peirce $1$-space is $\iota_1(S)\oplus \iota_2(S)$. Denote by $T_{2S}$ this orthogonal triple system.
Then, in case $S=S_8$, $T_{2S_8}$ is an orthogonal triple system defined on the direct sum of two copies of the split octonions, and we obtain:

\begin{theorem}\label{th:T2S}
Let $k$ be a field of characteristic $3$. Then the Lie superalgebra $\frg(T_{2S_8})$ of the orthogonal triple system $T_{2S_8}$ is isomorphic to $\frel(5;3)$.
\end{theorem}
\begin{proof}
Let $C$ be the split Cayley algebra over $k$, whose associated para-Hurwitz algebra is $S_8$, and let $J$ be the degree three simple Jordan algebra $H_3(C,*)$ considered above. Then, as vector spaces, $T_{2S_8}$ coincides with the Peirce $1$-space $J_1(e_0)$. The decomposition in \eqref{eq:H3Z2Z2} is a grading over $\bZ_2\times\bZ_2$ of the Jordan algebra $J$, and thus the Lie algebra of derivations of $J$ is also $\bZ_2\times\bZ_2$-graded as follows (see \cite[(3.12)]{CE2}):
\[
\der J=D_{\tri(S_8)}\oplus\bigl(\oplus_{i=0}^2D_i(S_8)\bigr),
\]
where, for $(d_0,d_1,d_2)\in \tri(S_8)$,
\[
\left\{\begin{aligned} &D_{(d_0,d_1,d_2)}(e_i)=0,\\
   &D_{(d_0,d_1,d_2)}\bigl(\iota_i(a)\bigr)=
   \iota_i\bigl(d_i(a)\bigr)
   \end{aligned}\right.
\]
for any $i=0,1,2$ and $a\in S_8$ (see \cite[(3.6)]{CE2}), while
\[
D_i(a)=2\bigl[ L_{\iota_i(a)},L_{e_{i+1}}\bigr]
\]
(indices modulo $3$) for $0\leq i\leq 2$ and $a\in S_8$, where $L_x$ denotes the left multiplication by $x$.

Note that $D_{\tri(S_8)}\oplus D_0(S_8)$ leaves $J_1(e_0)=\iota_1(S_8)\oplus\iota_2(S_8)$ invariant, and therefore embeds naturally in $\der T_{2S_8}$.

Besides, the Lie superalgebra of the orthogonal triple system $\hat J$ is (see \cite[\S 4]{CE2}):
\[
\frg(J)=\bigl(\frsp(V)\oplus\der J\bigr)\oplus\bigl(V\otimes \hat J),
\]
which is shown in \cite[Theorem 4.9]{CE2} to be isomorphic to $\frg(S_8,S_{1,2})$. Under this isomorphism $V\otimes T_{2S_8}$ corresponds to $\iota_1(S_8\otimes V)\oplus\iota_2(S_8\otimes V)$ inside $\frg(S_8,S_{1,2})$ which, under the isomorphism in Theorem \ref{th:el53maximalg83}, corresponds to the odd part of $\frel(5;3)$, and this odd part generates $\frel(5;3)$ as a Lie superalgebra. Therefore, the Lie superalgebra generated by $V\otimes T_{2S_8}$ corresponds to the subalgebra $\frg(S_8,S_{1,2})_+$ (isomorphic to $\frel(5;3)$). Using the isomorphism in \cite[Theorem 4.9]{CE2}, this proves that the subalgebra generated by $V\otimes T_{2S_8}$ in $\frg(J)$ is
\[
\bigl(\frsp(V)\oplus(D_{\tri(S_8)}\oplus D_0(S_8))\bigr)\oplus\bigl(V\otimes T_{2S_8}\bigr).
\]
Since this is a simple Lie superalgebra, by  Proposition \ref{pr:propertiesOTS}\,(2) it follows that it is isomorphic to the Lie superalgebra of the \OTS\ $T_{2S_8}$.
\end{proof}

\smallskip

\begin{remark}
Proposition \ref{pr:propertiesOTS}\,(3) shows that $\frtg(T_{2S_8})$ is a simple Lie algebra. By the proof above, it is $\bZ_2$-graded with even component isomorphic to $D_{\tri(S_8)}\oplus D_0(S_8)$, which is isomorphic to the orthogonal Lie algebra $\frso_9$, and with odd component (in the $\bZ_2$-grading) given by $T_{2S_8}$, which is the spin module for the even component. It follows that $\frtg(T_{2S_8})$ is the exceptional Lie algebra of type $F_4$.
\end{remark}

\bigskip

\section{Orthosymplectic triple systems and the Lie superalgebra $\frbr(2;3)$}

Orthosymplectic triple systems are the superversion of the orthogonal triple systems. They unify both orthogonal and symplectic triple systems. The definition was given in \cite[Definition 6.2]{CE2}:

\begin{definition}
Let $T=T\subo\oplus T\subuno$ be a vector superspace endowed with an even nonzero
supersymmetric bilinear form $(.\vert.):T\times T\rightarrow k$
(that is, $(T\subo\vert T\subuno)=0$, $(.\vert .)$ is symmetric on
$T\subo$ and alternating on $T\subuno$) and a triple product
$T\times T\times T\rightarrow T$: $(x,y,z)\mapsto [xyz]$
($[x_iy_jz_k]\in T_{i+j+k}$ for any $x_i\in T_i$, $y_j\in T_j$,
$z\in T_k$, where $i,j,k=\bar 0$ or $\bar 1$). Then $T$ is said to be an
\emph{orthosymplectic triple system} if it satisfies the following
identities:
\begin{subequations}\label{eq:OSTS}
\begin{align}
&[xyz]+(-1)^{  x   y }[yxz]=0\label{eq:OSTSa}\\
&[xyz]+(-1)^{  y   z }[xzy]=
 (x\vert y)z+(-1)^{  y   z }(x\vert z)y
 -2(y\vert z)x\label{eq:OSTSb}\\
&[xy[uvw]]=[[xyu]vw]+
 (-1)^{(  x +  y )  u }[u[xyv]w]+
 (-1)^{(  x +  y )(  u
    +  v )}[uv[xyw]]\label{eq:OSTSc}\\
&([xyu]\vert v)+(-1)^{(  x +  y )
         u }(u\vert [xyv])=0\label{eq:OSTSd}
\end{align}
\end{subequations}
for any homogeneous elements $x,y,u,v,w\in T$.
\end{definition}

\begin{remark}
If $T\subuno=0$, this is just the definition of an orthogonal triple system, while if $T\subo=0$, then it reduces to a symplectic triple system.
\end{remark}

\smallskip

As for orthogonal triple systems, the subspace $\inder T=\espan{[xy.]: x,y\in
T}$ is a subalgebra (actually an ideal) of the Lie superalgebra $\der T$
of derivations of $T$, whose elements are called \emph{inner
derivations}.

\begin{proposition}\label{pr:OSTS}
Let $T$ be a simple orthosymplectic triple system. Then its supersymmetric bilinear form $(.\vert .)$ is nondegenerate.
The converse is valid unless the characteristic of $k$ is $3$, $T=T\subuno$ and $\dim T=2$.
\end{proposition}
\begin{proof}
Given an orthosymplectic triple system, the kernel of its supersymmetric bilinear form: $T^\perp=\{ x\in T: (x\vert T)=0\}$, satisfies $[TTT^\perp]\subseteq T^\perp$ because of \eqref{eq:OSTSd}, while equations \eqref{eq:OSTSa} and \eqref{eq:OSTSb} show that $[TT^\perp T]=[T^\perp TT]\subseteq [TTT^\perp]+T^\perp$, so $T^\perp$ is an ideal of $T$. Thus, if $T$ is simple, then $(.\vert .)$ is nondegenerate.

Conversely, assume $T^\perp=0$ and let $I=I\subo\oplus I\subuno$ be a proper ideal of $T$. For homogeneous elements $x,y,z\in T$, \eqref{eq:OSTSb} shows that the element
\[
(x\vert y)z+(-1)^{yz}(x\vert z)y-2(y\vert z)x
\]
belongs to $I$ if at least one of $x,y,z$ is in $I$. For $x\in I$ we obtain
\[
(x\vert y)z+(-1)^{yz}(x\vert z)y\in I,
\]
while for $y\in I$, after permuting $x$ and $y$,
\[
(x\vert y)z-2(-1)^{yz}(x\vert z)y\in I,
\]
for homogeneous $x\in I$, $y,z\in T$. If the characteristic of $k$ is not $3$, it follows that $(I\vert T)T\subseteq I$, so $I=T$, a contradiction. But, even if the characteristic is $3$, it follows that the codimension $1$ subspace $(kx)^\perp=\{y\in T: (x\vert y)=0\}$ is contained in $I$ for any homogeneous element $x\in I$, and $I=T$ unless $\dim T=2$. In the latter case, either $T=T\subo$ or $T=T\subuno$. But for $T=T\subo$, $(x\vert y)y\in I$ for any homogeneous $x\in I$ and $y\in T$, and hence also $\{y\in T: (x\vert y)\ne 0\}$ is contained in $I$, so $I=T$. Thus $T=T\subuno$.
\end{proof}

\begin{remark} The two dimensional symplectic triple system in \cite[Proposition 2.7(i)]{EldNew} shows that there are indeed nonsimple orthosymplectic triple systems of superdimension $0\vert 2$ (that is, $\dim T\subo=0$, $\dim T\subuno=2$).
\end{remark}

As for orthogonal triple systems, the following result (see \cite[Theorem 6.3]{CE2}) holds:

\begin{proposition}\label{pr:gOSTS}
Let $\bigl(T,[...],(.\vert.)\bigr)$ be an orthosymplectic triple system and let
$\bigl(V,\langle.\vert.\rangle\bigr)$ be a two dimensional vector
space endowed with a nonzero alternating bilinear form. Let $\frs$
be a Lie subsuperalgebra of $\der T$ containing $\inder T$. Define
the $\bZ_2$-graded superalgebra
$\frg=\frg(T,\frs)=\frg(0)\oplus\frg(1)$ with
\[
\begin{cases}
\frg(0)= \frsp(V)\oplus \frs&\text{(so $\frg(0)\subo=\frsp(V)\oplus\frs\subo$,
$\frg(0)\subuno=\frs\subuno$),}\\
\frg(1)=V\otimes T&\text{(with $\frg(1)\subo=V\otimes T\subuno$,
$\frg(1)\subuno=V\otimes T\subo$, $V$ is odd!)},
\end{cases}
\]
and superanticommutative multiplication given by:
\begin{itemize}
\item
the multiplication on $\frg(0)$ coincides with its bracket as a Lie superalgebra;
\item
$\frg(0)$ acts naturally on $\frg(1)$:
\[
[s,v\otimes x]=s(v)\otimes x,\qquad [d,v\otimes x]=(-1)^{
d }v\otimes d(x),
\]
for any $s\in \frsp(V)$, $v\in V$, and homogeneous elements $d\in\frs$ and $x\in T$;
\item
for any $u,v\in V$ and homogeneous $x,y\in T$:
\begin{equation}\label{eq:gTosproduct}
[u\otimes x,v\otimes y]=(-1)^{  x }\bigl(
  -(x\vert y)\gamma_{u,v} +\langle u\vert v\rangle d_{x,y}\bigr)
\end{equation}
where $\gamma_{u,v}=\langle u\vert .\rangle v+\langle v\vert
.\rangle u$ and $d_{x,y}=[xy.]$.
\end{itemize}
Then $\frg(T,\frs)$ is a $\bZ_2$-graded Lie superalgebra. Moreover,
$\frg(T,\frs)$ is simple if and only if $\frs$ coincides with
$\inder T$ and $(.\vert.)$ is nondegenerate.

Conversely, given a $\bZ_2$-graded Lie superalgebra
$\frg=\frg(0)\oplus\frg(1)$ with
\[
\begin{cases}
\frg(0)=\frsp(V)\oplus \frs,\\
\frg(1)=V\otimes T,
\end{cases}
\]
where $T$ is an $\frs$-module and $V$ is
considered as an odd vector space, by $\frsp(V)$-invariance of the
 bracket, equation \eqref{eq:gTosproduct} is satisfied for an even
supersymmetric bilinear form $(.\vert .):T\times T\rightarrow k$ and
a superantisymmetric bilinear map $d_{.,.}:T\times T\rightarrow
\frs$. Then, if $(.\vert .)$ is not $0$ and a triple product on $T$
is defined by means of $[xyz]=d_{x,y}(z)$, $T$ becomes and
orthosymplectic triple system and the image of $\frs$ in $\frgl(T)$
under the given representation is a subalgebra of $\der T$
containing $\inder T$.
\end{proposition}

The Lie superalgebra $\frg(T)=\frg(T,\inder(T))$ is called the \emph{Lie superalgebra of the orthosymplectic triple system $T$}.

If the characteristic of the ground field $k$ is equal to $3$, then for any homogeneous elements $x,y,z$ in an orthosymplectic triple system, we have:
\[
\begin{split}
[xyz]+&(-1)^{x(y+z)}[yzx]+(-1)^{(x+y)z}[zxy]\\
 &=[xyz]+(-1)^{x(y+z)}[yzx]-2(-1)^{(x+y)z}[zxy]\\
 &=\bigl([xyz]+(-1)^{yz}[xzy]\bigr)\\
 &\qquad\qquad -(-1)^{xy+xz+yz}\bigl([zyx]+(-1)^{xy}[zxy]\bigr)\quad\textrm{(by \eqref{eq:OSTSa})}\\
 &=\bigl((x\vert y)z+(-1)^{yz}(x\vert z)y-2(y\vert z)x\bigr)\\
 &\qquad\qquad-(-1)^{xy+xz+yz}\bigl((z\vert y)+(-1)^{xy}(z\vert x)y-2(y\vert x)z)\quad\textrm{(by \eqref{eq:OSTSb})}\\
 &=3\bigl((x\vert y)z-(y\vert z)z\bigr)=0,
\end{split}
\]
so that, as in \cite[Theorem 5.1]{EldNew}:

\begin{proposition}\label{pr:gtOSTS}
Let $\bigl(T,[...],(.\vert .)\bigr)$ be an orthosymplectic triple system over a field $k$ of characteristic $3$. Define the $\bZ_2$-graded superalgebra $\frtg=\frtg(T)=\frtg(T)_+\oplus\frtg(T)_-$, with
\[
\frtg_+=\inder(T),\qquad\frtg_-=T,
\]
and superanticommutative multiplication given by:
\begin{itemize}
\item
the multiplication on $\frtg_+$ coincides with its bracket as a Lie superalgebra;
\item
$\frtg_+$ acts naturally on $\frtg_-$, that is,
$[d,x]=d(x)$ for any $d\in\inder(T)$ and $x\in T$;
\item
$[x,y]=d_{x,y}=[xy.]$, for any $x,y\in \frtg_-=T$.
\end{itemize}
Then $\frtg$ is a $\bZ_2$-graded Lie superalgebra, with the even part $\frtg\subo=\inder(T)\subo\oplus T\subo$ and the odd part $\frtg\subuno=\inder(T)\subuno\oplus T\subuno$. Moreover, $T$ is a simple orthosymplectic triple system if and only if  $\frtg$ is simple as a $\bZ_2$-graded Lie superalgebra.
\end{proposition}

\medskip

Now, let $C$ be a Hurwitz superalgebra of dimension $>1$ over a field $k$ of characteristic $\ne 2$, with norm $q=(q\subo,\bil)$, and standard involution $x\mapsto \bar x$. For any homogeneous elements $x,y,z$, the following holds (see \cite{EldOkuCompoSuper}):
\[
\begin{split}
&\bil(xy,z)=(-1)^{xy}\bil(y,\bar xz)=(-1)^{yz}\bil(x,z\bar y),\\
&x\bar y+(-1)^{xy}y\bar x=\bil(x,y)1=\bar xy+(-1)^{xy}\bar yx,\\
&\bar x(yz)+(-1)^{xy}\bar y(xz)=\bil(x,y)z=(zx)\bar y+(-1)^{xy}(zy)\bar x.
\end{split}
\]
Consider the subspace of trace zero elements,
\[
C^0=\{ x\in C:\bil(1,x)=0\}=\{x\in C: \bar x=-x\}.
\]
Then, for any homogeneous elements $x,y\in C^0$, we have
\[
xy+(-1)^{xy}yx=-(x\bar y+(-1)^{xy}y\bar x)=-\bil(x,y)1,
\]
while $xy-(-1)^{xy}yx=[x,y]$. Thus
\begin{equation}\label{eq:xyb[xy]}
xy=\frac{1}{2}\bigl(-\bil(x,y)1+[x,y]\bigr).
\end{equation}
Also, for any homogeneous elements $x,y,z\in C^0$, we have
\[
\begin{split}
\bil([x,y],z)&=\bil(xy-(-1)^{xy}yx,z)\\
 &=\bil(x,(-1)^{yz}z\bar y-\bar yz)\\
 &=\bil(x,yz-(-1)^{yz}zy)\\
 &=\bil(x,[y,z]),
\end{split}
\]
so
\begin{equation}\label{eq:bilassociative}
\bil([x,y],z)=\bil(x,[y,z])
\end{equation}
for any $x,y,z\in C^0$. Using \eqref{eq:xyb[xy]} and \eqref{eq:bilassociative} we obtain:
\[
\begin{split}
[[x,y],z]+&(-1)^{yz}[[x,z],y]\\
 &=\bil([x,y],z)1+2[x,y]z+(-1)^{yz}\bigl(\bil([x,z],y)1+2[x,z]y\bigr)\\
 &=2\bigl([x,y]z+(-1)^{yz}[x,z]y\bigr)\\
 &=2\bigl(\bil(x,y)z+2(xy)z+(-1)^{yz}\bigl(\bil(x,z)y+2(xz)y\bigr)\bigr)\\
 &=2\bigl(\bil(x,y)z+(-1)^{yz}\bil(x,z)y\bigr)-4\bigl((xy)\bar z+(-1)^{yz}(xz)\bar y\bigr)\\
 &=2\bil(x,y)z+2(-1)^{yz}\bil(x,z)y-4\bil(y,z)x,
\end{split}
\]
for any homogeneous $x,y,z\in C^0$. Therefore, with $(x\vert y)=2\bil(x,y)$, for any $x,y,z\in C^0$ we have:
\begin{equation}\label{eq:C0OSTS}
[[x,y],z]+(-1)^{yz}[[x,z],y]=(x\vert y)z+(-1)^{yz}(x\vert z)y-2(y\vert z)x.
\end{equation}

Now, if the characteristic of the ground field $k$ is equal to $3$, for any homogeneous $x,y,z\in C^0$ we have:
\[
\begin{split}
[[x,y],z]+&(-1)^{x(y+z)}[[y,z],x]+(-1)^{(x+y)z}[[z,x],y]\\
 &=[[x,y],z]+(-1)^{x(y+z)}[[y,z],x]-2(-1)^{(x+y)z}[[z,x],y]\\
 &=\bigl([[x,y],z]+(-1)^{yz}[[x,z],y]\bigr)\\
 &\qquad\qquad -(-1)^{xy+xz+yz}\bigl([[z,y],x]+(-1)^{xy}[[z,x],y]\bigr)\\
 &=\bigl((x\vert y)z+(-1)^{yz}(x\vert z)y-2(y\vert z)x\bigr)\\
 &\qquad\qquad-(-1)^{xy+xz+yz}\bigl((z\vert y)+(-1)^{xy}(z\vert x)y-2(y\vert x)z)\\
 &=3\bigl((x\vert y)z-(y\vert z)z\bigr)=0.
\end{split}
\]
Thus, $(C^0,[.,.])$ is a Lie superalgebra, and then equations \eqref{eq:bilassociative} and \eqref{eq:C0OSTS}
show the validity of the first assertion in the following result:

\begin{theorem}\label{th:C0OSTS}
Let $C$ be a Hurwitz superalgebra of dimension $\geq 2$ over a field $k$ of characteristic $3$ with norm $q=(q\subo,\bil)$. Then, with the triple product $[xyz]=[[x,y],z]$ and the supersymmetric bilinear form $(.\vert .)=2\bil(.,.)$, $C^0$ becomes an orthosymplectic triple system. Moreover, if the dimension of $C$ is $\leq 3$, then the triple product is trivial, otherwise the inner derivation algebra $\inder(C^0)$ equals $\ad_{C^0}$, the linear span of the adjoint maps $\ad_x\,(:y\mapsto [x,y])$ for any $x\in C^0$.
\end{theorem}
\begin{proof}
Only the last assertion needs to be verified. If the dimension of $C$ is at most $3$, then $C$ is supercommutative, so $[C^0,C^0]=0$. However, if the dimension of $C$ is at least $4$ (hence either $4$, $6$ or $8$) then $[C^0,C^0]=C^0$.
\end{proof}

\begin{corollary}\label{co:C0OSTS}
Let $C$ be a Hurwitz superalgebra of dimension $\geq 4$ over a field $k$ of characteristic $3$. Let $V$ be a two-dimensional vector space endowed with a nonzero alternating bilinear form $\langle.\vert .\rangle$. Consider the anticommutative superalgebra
\[
\frg=\bigl(\frsp(V)\oplus C^0)\oplus\bigl(V\otimes C^0\bigr),
\]
with $\frg\subo=(\bigl(\frsp(V)\oplus (C^0)\subo\bigr)\oplus\bigl(V\otimes (C^0)\subuno\bigr)$ and $\frg\subuno=(C^0)\subuno\oplus \bigl(V\otimes(C^0)\subo\bigr)$, and multiplication given by:
\begin{itemize}
\item the usual Lie bracket in the direct sum of the Lie algebra $\frsp(V)$ and the Lie superalgebra $C^0$,
\item $[\gamma,v\otimes x]=\gamma(v)\otimes x$, for any $\gamma\in \frsp(V)$, $v\in V$ and $x\in C^0$,
\item $[x,v\otimes y]=(-1)^xv\otimes [x,y]$, for any homogeneous $x,y\in C^0$ and $v\in V$,
\item $[u\otimes x, v\otimes y]=(-1)^x\bigl(-(x\vert y)\gamma_{u,v}+\langle u\vert v\rangle [x,y]\bigr)\, \in\frsp(V)\oplus C^0$, for any $u,v\in V$ and homogeneous $x,y\in C^0$ (where, as before, $(x\vert y)=2\bil(x,y)$ and $\gamma_{u,v}=\langle u\vert .\rangle v+\langle v\vert .\rangle u$).
\end{itemize}
Then $\frg$ is a Lie superalgebra.
\end{corollary}
\begin{proof}
It suffices to note that the Lie superalgebra $\frg$ is just the Lie superalgebra $\frg(C^0,\inder(C^0))$ in Proposition \ref{pr:gOSTS} of the orthosymplectic triple system $\bigl(C^0,[...],(.\vert.)\bigr)$ after the natural identification of $\inder(C^0)=\ad_{C^0}$ with $C^0$.
\end{proof}

\medskip

If the dimension of $C$ in Corollary \ref{co:C0OSTS} is $4$ (and hence $C$ is a quaternion algebra), it is not difficult to see that the Lie superalgebra $\frg$ is a form of the orthosymplectic Lie superalgebra $\frosp_{3,2}$. Also, if the dimension of $C$ is $8$, so that $C$ is an algebra of octonions, then $\frg$ is a form of the Lie superalgebra that appears in \cite[Theorem 4.22(i)]{EldNew}, which is the derived subalgebra of the Lie superalgebra $\frg(S_2,S_{1,2})$ in the Supermagic Square (see \cite[Corollary 4.10(ii)]{CE2} and \cite[\S 3]{EldTits3}). Also note that if the characteristic is not $3$, then $C^0$ is still an orthogonal triple system, but its associated Lie superalgebra is a simple Lie superalgebra of type $G(3)$ (see \cite[Theorem 4.7 (G-type)]{EldNew}).

\smallskip

We are left with the $4\vert 2$ dimensional Hurwitz superalgebra $C=B(4,2)$ in \eqref{eq:B42} over a field $k$ of characteristic $3$. The Lie bracket of elements in $C^0$ is given by:
\begin{itemize}
\item The usual bracket $[f,g]=fg-gf$ in $\frsl(V)=\frsp(V)$.

\item $[f,u]=f\cdot u-u\cdot f=-2f(u)=f(u)$ for any $f\in \frsp(V)$ and $u\in V$.

\item $[u,v]=u\cdot v-(-1)^{uv}v\cdot u=u\cdot v+v\cdot u=\bil(.,u)v+\bil(.,v)u=(u\vert .)v+(v\vert .)u$ for any $u,v\in V$ (recall that $(.\vert .)=2\bil(.,.)=-\bil(.,.)$).
\end{itemize}

\begin{proposition}
The Lie superalgebra $B(4,2)^0$ is isomorphic to the orthosymplectic Lie superalgebra $\frosp_{1,2}$.
\end{proposition}
\begin{proof}
The orthosymplectic Lie superalgebra $\frosp_{1,2}$ is the subalgebra of the general Lie superalgebra $\frgl(1,2)$ given by:
\[
\frosp_{1,2}=\left\{\left(\begin{array}{c|cc} 0&-\nu&\mu\\
\hline \mu&\alpha&\beta\\ \nu&\gamma&-\alpha\end{array}\right): \alpha,\beta,\gamma,\mu,\nu\in k\right\}.
\]
Fix a basis $\{u,v\}$ of $V$ with $(u\vert v)=1$, and consider the linear map:
\[
\begin{split}
C^0=\frsp(V)\oplus V& \longrightarrow \frosp_{1,2}\\
f\in \frsp(V)&\mapsto \left(\begin{array}{c|cc} 0&0&0\\
\hline 0&\alpha&\beta\\ 0&\gamma&-\alpha\end{array}\right)\quad\textrm{with}\quad\begin{cases} f(u)=\alpha u+\gamma v,&\\ f(v)=\beta u -\alpha v,&\end{cases}\\
\mu u+\nu v\in V&\mapsto \left(\begin{array}{c|cc} 0&-\nu&\mu\\
\hline \mu&0&0\\ \nu&0&0\end{array}\right)
\end{split}
\]
This is checked to be an isomorphism of Lie algebras by straightforward computations.
\end{proof}

Also note that for $f\in \frsp(V)$, $f^2=-\det(f)1$ (by the Cayley-Hamilton equation) and $q\subo(f)=\det(f)$ and $\tr(f^2)=-2\det(f)=\det(f)$, so $q\subo(f)=\tr(f^2)$,  $\bil(f,g)=2\tr(fg)$, and $(f\vert g)=\tr(fg)$ for any $f,g\in \frsp(V)=(C^0)\subo$. (Here $\tr$ denotes the usual trace in $\End_k(V)=B(4,2)\subo$).

\begin{theorem}\label{th:OSTSbr23}
The Lie superalgebra of the orthosymplectic triple system $B(4,2)^0$ is isomorphic to the Lie superalgebra $\frbr(2;3)$.
\end{theorem}
\begin{proof}
Since there are two vector spaces of dimension $2$ involved here, let us denote them by $V_1$ and $V_2$, whose nonzero alternating bilinear forms will be both denoted by $\langle.\vert.\rangle$. Then consider the Hurwitz superalgebra $C=B(4,2)=\End_k(V_2)\oplus V_2$, as defined in \eqref{eq:B42}. The Lie superalgebra associated to the orthosymplectic triple system $C^0$ is given, up to isomorphism, in Corollary \ref{co:C0OSTS}:
\[
\frg=\bigl(\frsp(V_1)\oplus C^0\big)\oplus \bigl(V_1\otimes C^0\bigr).
\]
Its even part is
\[
\frg\subo=\bigl(\frsp(V_1)\oplus\frsp(V_2)\bigr)\oplus \bigl(V_1\otimes V_2),
\]
with multiplication given by the natural Lie bracket in the direct sum $\frsp(V_1)\oplus\frsp(V_2)$, the natural action of this subalgebra on $V_1\otimes V_2$, and by
\[
[a\otimes u,b\otimes v]=(u\vert v)\gamma_{a,b}-\langle a\vert b\rangle\gamma_{u,v},
\]
for any $a,b\in V_1$ and $u,v\in V_2$, where $(.\vert .)=2\bil(.,.)$. Here $\gamma_{a,b}=\langle a\vert .\rangle b +\langle b\vert .\rangle a$, while $\gamma_{u,v}=(u\vert .)v+(v\vert .)u$. This Lie algebra is precisely the Lie algebra $L(1)$ of Kostrikin (see \cite{Kostrikin} or \cite[Proposition 2.12]{EldNew}).

On the other hand, its odd part is
\[
\frg\subuno=V_2\oplus\bigl(V_1\otimes \frsp(V_2)\bigr).
\]
Since $C^0$ is a simple orthosymplectic triple system, the Lie superalgebra $\frg$ is simple (Proposition \ref{pr:gOSTS}). Fix bases $\{a_i,b_i\}$ of $V_i$ ($i=1,2$) with $\langle a_1\vert b_1\rangle=1=(a_2\vert b_2)$, and let $h_i,e_i,f_i\in \frsp(V_i)$ be given by
\begin{equation}\label{eq:hef}
\begin{aligned}
h_i(a_i)&=a_i,&\qquad h_i(b_i)&=-b_i,\\
e_i(a_i)&=0,& e_i(b_i)&=a_i,\\
f_i(a_i)&=b_i,& f_i(b_i)&=0.
\end{aligned}
\end{equation}
Then $\espan{h_1,h_2}$ is a Cartan subalgebra of $\frg$, and it is the $(0,0)$-component of the $\bZ\times\bZ$-grading of $\frg$ obtained by assigning $\deg(a_i)=\epsilon_i$, $\deg(b_i)=-\epsilon_i$, $i=1,2$, where $\{\epsilon_1,\epsilon_2\}$ is the canonical basis of $\bZ\times\bZ$. The set of nonzero degrees is
\[
\Phi=\{\pm 2\epsilon_1,\pm 2\epsilon_2,\pm\epsilon_1\pm\epsilon_2,\pm\epsilon_2,\pm\epsilon_1, \pm\epsilon_1\pm2\epsilon_2\}.
\]
Consider the elements
\begin{alignat*}{3}
E_1&=a_1\otimes f_2,&\quad F_1&=b_1\otimes e_2,&\quad H_1&=[E_1,F_1]=h_1-h_2,\\
E_2&=a_2,& F_2&=-b_2,& H_2&=[E_2,F_2]=h_2.
\end{alignat*}
Then we have that the subspace $\espan{H_1,H_2}=\espan{h_1,h_2}$ is the previous Cartan subalgebra of $\frg$, $E_1$ belongs to the homogeneous component $\frg_{\epsilon_1-2\epsilon_2}$ in the $\bZ\times\bZ$-grading, and similarly $F_1\in\frg_{-\epsilon_1+2\epsilon_2}$, $E_2\in \frg_{\epsilon_2}$, and $F_2\in\frg_{-\epsilon_2}$. Also, the elements $E_1,E_2,F_1,F_2$ generate the Lie superalgebra $\frg$. Besides,
\[
\begin{split}
[H_1,E_1]&=h_1(a_1)\otimes f_2-a_1\otimes [h_2,f_2]=a_1\otimes f_2+2a_1\otimes f_2=0,\\
[H_1,E_2]&=(h_1-h_2)(a_2)=-a_2,\\
[H_2,E_1]&=a_1\otimes [h_2,f_2]=-2a_1\otimes f_2,\\
[H_2,E_2]&=h_2(a_2)=a_2,
\end{split}
\]
and similarly for the action of the $H_i$'s on the $F_j$'s. It follows, with the same arguments as in \cite[\S 4]{CE1}, that $\frg$ is the Lie superalgebra with Cartan matrix $\left(\begin{smallmatrix} 0&-1\\ -2&1\end{smallmatrix}\right)$, which is the first Cartan matrix of the Lie superalgebra $\frbr(2;3)$ given in \cite[\S 10.1]{BGLCartan}.
\end{proof}

In this way, the Lie superalgebra $\frbr(2;3)$, of superdimension $10|8$ is completely determined by the orthosymplectic triple system $B(4,2)^0$ (that is, by the orthosymplectic triple system obtained naturally from the Lie superalgebra $\frosp_{1,2}$) of superdimension $3|2$.

\bigskip

\section{Orthosymplectic triple systems and the Lie superalgebra $\frel(5;3)$}

In this section, the characteristic of the ground field $k$ will always be assumed to be $3$, since we will be dealing with the superalgebras $S_{1,2}$ and $\frel(5;3)$, which only make sense in this characteristic.

\smallskip

Equation \eqref{eq:g8301}, together with Theorem \ref{th:el53maximalg83}, which allows us to identify the Lie superalgebra $\frel(5;3)$ with the maximal subalgebra $\frg(S_8,S_{1,2})_+$, show that there is a decomposition of $\frg(S_8,S_{1,2})$ into the direct sum ($\bZ_2$-grading) $\frg(S_8,S_{1,2})=\frel(5;3)\oplus T$, where:
\begin{equation}\label{eq:el53T}
\begin{split}
\frel(5;3)&=\bigl(\tri(S_8)\oplus\frsp(V)\bigr) \oplus \iota_0(S_8\otimes 1)\oplus \bigl(\iota_1(S_8\otimes V)\oplus \iota_2(S_8\otimes V)\bigr),\\
T&=\bigl(\iota_1(S_8\otimes 1)\oplus\iota_2(S_8\otimes 1)\bigr)\oplus\bigl(V\oplus \iota_0(S_8\otimes V)\bigr),\\
\end{split}
\end{equation}
where $V$ is a two dimensional vector space endowed with a nonzero alternating bilinear form.

This section will show that $T$ is an orthosymplectic triple system, with the triple product given by $[xyz]=[[x,y],z]$ and a suitable supersymmetric bilinear form, and that $\frel(5;3)$ is isomorphic to the Lie superalgebra of derivations of this orthosymplectic triple system.

A few preliminary results are needed.

\begin{lemma}\label{le:Bilg83}
There exists a unique supersymmetric associative bilinear form
\[
B:\frg(S_8,S_{1,2})\times \frg(S_8,S_{1,2})\rightarrow k
\]
such that
\begin{equation}\label{eq:Biotas}
B\bigl(\iota_i(x\otimes u),\iota_j(y\otimes v)\bigr) =\delta_{ij}\bil(x,y)\bil(u,v),
\end{equation}
for any $i,j=0,1,2$,  $x,y\in S_8$ and $u,v\in S_{1,2}$. (Here $\delta_{ij}$ is the usual Kronecker delta and $\bil$ denotes the polar form of the norm in both $S_8$ and $S_{1,2}$.)
\end{lemma}
\begin{proof}
This is proved as in \cite[Corollary 4.9]{NewLook}. First there is a unique invariant supersymmetric bilinear form $B_{1,2}$ on the orthosymplectic Lie superalgebra $\frosp(S_{1,2},q)$ such that:
\[
B_{1,2}(d,\sigma_{u,v})=\bil(d(u),v)
\]
for any $u,v\in S_{1,2}$ and $d\in \frosp(S_{1,2},q)$, where $\sigma_{u,v}$ is defined in \eqref{eq:sigmaxy}. Actually, $B_{1,2}$ is given by $B_{1,2}(d,d')=-\frac{1}{2}\str(dd')$, where $\str$ denotes the supertrace. (Note that
\[
\str(\sigma_{x,y}\sigma_{u,v}) =-2\bigl((-1)^{yu}\bil(x,u)\bil(y,v)-(-1)^{(y+u)v}\bil(x,v)\bil(y,u)\bigr).)
\]
Also, in \cite{NewLook} it is proved that there is a unique invariant symmetric bilinear form $B_8$ on $\tri(S_8)$ such that
\[
B_8\bigl((d_0,d_1,d_2),\theta^i(t_{x,y})\bigr)=\bil(d_i(x),y)
\]
for any $x,y\in S_8$ and $(d_0,d_1,d_2)\in \tri(S_8)$.

Then the supersymmetric invariant bilinear form $B$ required is defined by imposing the following conditions:
\begin{itemize}
\item The restriction of $B$ to $\tri(S_{1,2})$ is given by $B_{1,2}$ (after identifying $\tri(S_{1,2})$ with $\frosp(S_{1,2},q)$ because of \eqref{eq:triS12}).
\item The restriction of $B$ to $\tri(S_8)$ is given by $B_8$.
\item The restriction of $B$ to $\oplus_{i=0}^2\iota_i(S_8\otimes S_{1,2})$ is given by \eqref{eq:Biotas}.
\end{itemize}
\end{proof}

Note that $\frg(S_8,S_{1,2})$ is then the orthogonal direct sum, relative to $B$, of the subspaces $\tri(S_8)$, $\tri(S_{1,2})$ and $\iota_i(S_8\otimes S_{1,2})$, $i=0,1,2$.

Now, the description of $\frg(S_8,S_{1,2})$ in the proof of Theorem \ref{th:el53maximalg83} becomes quite useful in the proof of the next result:

\begin{lemma}\label{le:derg83}
Any derivation of the Lie superalgebra $\frg(S_8,S_{1,2})$ is inner.
\end{lemma}
\begin{proof}
As in \cite{CE1}, take five two-dimensional vector spaces $V_1,\ldots,V_5$ endowed with nonzero alternating bilinear forms $\langle.\vert .\rangle$. Take symplectic bases $\{v_i,w_i\}$ of $V_i$ for any $i=1,\ldots,5$ ($\langle v_i\vert w_i\rangle =1$) and the basis $\{h_i,e_i,f_i\}$ of $\frsp(V_i)$ given by
\[
h_i=\gamma_{v_i,w_i},\quad e_i=\gamma_{w_i,w_i},\quad f_i=-\gamma_{v_i,v_i},
\]
which satisfy $[h_i,e_i]=2e_i$, $[h_i,f_i]=-2f_i$, and $[e_i,f_i]=h_i$.

Consider the description of $\frg(S_8,S_{1,2})$ in \eqref{eq:g83sigmas}: $\frg(S_8,S_{1,2})=\oplus_{\sigma\in\calS_{8,3}}V(\sigma)$. This shows that $\frg(S_8,S_{1,2})$ is $\bZ^5$-graded, by assigning $\deg w_i=\epsilon_i$, $\deg v_i=-\epsilon_i$, where $\{\epsilon_1,\ldots,\epsilon_5\}$ is the canonical basis of $\bZ^5$. The vector subspace $\frh=\espan{h_1,\ldots,h_5}$ is a Cartan subalgebra of $\frg(S_8,S_{1,2})$. Consider the $\bZ$-linear map:
\[
\begin{split}
R:\bZ^5&\longrightarrow \frh^*\\
\epsilon_i&\mapsto\  R(\epsilon_i):h_j\mapsto \delta_{ij}.
\end{split}
\]
The set of nonzero degrees of $\frg(S_8,S_{1,2})$ in the $\bZ^5$-grading is given by
\begin{multline*}
\Phi=\bigl\{\pm 2\epsilon_i:i=1,\ldots,5\bigr\}\\
\cup \bigl\{\pm \epsilon_{i_1}\pm\cdots\pm\epsilon_{i_r}: 1\leq i_1<\cdots<i_r\leq 5,\ \{i_1,\ldots,i_r\}\in\calS_{8,3}\setminus\{\emptyset\}\bigr\}.
\end{multline*}
The set $R(\Phi)$ is the set of roots of $\frg(S_8,S_{1,2})$ relative to the Cartan subalgebra $\frh$. Note that the restriction of $R$ to $\Phi$ fails to be one-to-one only because $\{\pm 2\epsilon_5,\pm\epsilon_5\}$ is contained in $\Phi$, and $R(\pm 2\epsilon_5)=R(\mp\epsilon_5)$, as the characteristic is equal to $3$.

The Lie superalgebra of derivations of $\frg=\frg(S_8,S_{1,2})$ inherits the $\bZ^5$-grading, so in order to prove the Lemma it is enough to prove that homogeneous derivations (in this $\bZ^5$-grading) are inner. Thus, assume that $d\in \der(\frg)_\nu$, with $\nu\in\bZ^5$:
\begin{enumerate}
\item If $\nu\ne 0$ and $d(\frh)=0$ (note that the Cartan subalgebra $\frh$ is just the $0$-component in this grading), then $d$ preserves the eigenspaces (root spaces) of $\frh$, and hence $d(\frg_\mu)=0$ for any $\mu\in \Phi\setminus\{\pm 2\epsilon_5,\pm\epsilon_5\}$, as $d(\frg_\mu)$ must simultaneously be contained in $\frg_{\mu+\nu}$ and in the root space of root $R(\mu)$. But the subspaces $\frg_\mu$, with $\mu\in \Phi\setminus\{\pm 2\epsilon_5,\pm\epsilon_5\}$ generate the Lie superalgebra $\frg$. (This can be checked easily, but it also follows from \cite[Proposition 5.25]{CE1}.) Hence $d=0$, which is trivially inner.

\item If $\nu\ne 0$ and $d(\frh)\ne 0$, then $d(\frh)$ is contained in $\frg_\nu$, which has dimension at most $1$. Thus $\frg_\nu=kx_\nu$ for some $x_\nu$. Then for any $h\in \frh$, $d(h)=f(h)x_\nu$ for some $f\in \frh^*$. Then, for any $h,h'\in \frh$,
\[
0=d([h,h'])=[d(h),h']+[h,d(h')] =\bigl(-R(\nu)(h')f(h)+R(\nu(h)f(h')\bigr)x_\nu.
\]
    As $R(\nu)\ne 0$, it follows that there is a scalar $\alpha\in k$ with $f=\alpha R(\nu)$, and $d'=d-\alpha\ad x_\nu$ is another derivation in $\der(\frg)_\nu$ with $d'(\frh)=0$, so $d'$ must be $0$ by the previous case, and hence $d$ is inner.

\item Finally, if $\nu=0$, then $d(e_i)\in\frg_{2\epsilon_i}=ke_i$, so $d(e_i)=\alpha_ie_i$ for any $i$. Also, $d(f_i)=\beta_if_i$ for any $i$ ($\alpha_i,\beta_i\in k$). As $ke_i+kf_i+kh_i$ is a Lie subalgebra isomorphic to $\frsl_2$, it follows at once that $\alpha_i+\beta_i=0$. Then the derivation $d'=d-\frac{1}{2}\ad_{\alpha_1h_1+\cdots+\alpha_5h_5}$ satisfies $d'(e_i)=0=d'(f_i)$ for any $i$, so $d'(h_i)=0$, and hence $d'(\frsp(V_i))=0$ for any $i$. As $d'$ preserves degrees, it preserves each subspace $V(\sigma)$, for $\emptyset\ne\sigma\in \calS_{8,3}$, which is an irreducible module for $\oplus_{i=1}^5\frsp(V_i)$. By Schur's Lemma, there is a scalar $\alpha_\sigma\in k$ such that the restriction of $d'$ to any $V(\sigma)$ is $\alpha_\sigma id$. But $0\ne [V(\sigma),V(\sigma)]\subseteq \oplus_{i=1}^5\frsp(V_i)$, so $2\alpha_\sigma=0$ for any such $\sigma$ and $d'=0$. Thus $d$ is inner in this case, too.
\end{enumerate}
\end{proof}

\medskip

Consider now the triple product on the subspace $T$ (the odd component in the $\bZ_2$-grading of $\frg(S_8,S_{1,2})$ considered so far) inherited from the Lie bracket in $\frg(S_8,S_{1,2})$:
\[
\begin{split}
T\otimes T\otimes T&\rightarrow T\\
X\otimes Y\otimes Z&\mapsto [XYZ]=[[X,Y],Z],
\end{split}
\]
As $T$ is the odd component of $\frg(S_8,S_{1,2})$, it is a Lie triple supersystem. Therefore $(T,[...])$ satisfies equations \eqref{eq:OSTSa} and \eqref{eq:OSTSc}.

Also, if we consider the supersymmetric bilinear form $(.\vert.)$ on $T$ given by the restriction of the bilinear form $B$ given in Lemma \ref{le:Bilg83}, the invariance of $B$ immediately shows that $\bigl(T,[...],(.\vert.)\bigr)$ also satisfies equation \eqref{eq:OSTSd}.

\begin{theorem}
$\bigl(T,[...],(.\vert.)\bigr)$ is an orthosymplectic triple system whose Lie superalgebra of derivations is isomorphic to $\frel(5;3)$. Moreover, the associated Lie superalgebra $\frg(T)$ is isomorphic to the Lie superalgebra $\frg(S_{4,2},S_{4,2})$ in the Supermagic Square.
\end{theorem}
\begin{proof}
It is enough to check equation \eqref{eq:OSTSb}.

Take a symplectic basis $\{a,b\}$ of the two dimensional vector space $V$ in \eqref{eq:el53T} (that is, $\langle a\vert b\rangle=1$), then $T$ is generated, as a module over $\frel(5;3)$ by $\iota_0(S_8\otimes a)$ or by $\iota_0(S_8\otimes b)$. Also, $T\otimes T$ is generated by $\iota_0(S_8\otimes a)\otimes\iota_0(S_8\otimes b)$.
Both the left and right sides of equation \eqref{eq:OSTSb} are given by  $\frel(5;3)$-invariant trilinear maps $T\otimes T\otimes T\rightarrow T$. Therefore, it is enough to prove that the condition
\[
\begin{split}
[X\,\iota_0(y\otimes a)\,\iota_0(z\otimes b)]
 &-[X\,\iota_0(z\otimes b)\,\iota_0(y\otimes a)]\\
 &=\bigl(X\vert\iota(y\otimes a)\bigr)\iota_0(z\otimes b)-\bigl(X\vert\iota_0(z\otimes b)\bigr)\iota_0(y\otimes a)+\bil(y,z)X
\end{split}
\]
holds for any $X\in T$ and $y,z\in S_8$.

\begin{itemize}
\item For $X=u\in V\simeq\tri(S_{1,2})\subuno$, since $\{a,b\}$ is a symplectic basis, $u=\langle u\vert b\rangle a-\langle u\vert a\rangle b$, so:
\[
\begin{split}
[u\,\iota_0(y\otimes a)\,\iota_0(z\otimes b)]&=
  -\langle u\vert a\rangle[\iota_0(y\otimes 1),\iota_0(z\otimes b)]\\
  &=-\langle u\vert a\rangle\bil(y,z)t_{1,b}\\
  &=-\langle u\vert a\rangle \bil(y,z)b,
\end{split}
\]
where, as before, $V$ is identified with $\tri(S_{1,2})\subuno$. Thus,
\[
\begin{split}
[X\,\iota_0(y\otimes a)\,\iota_0(z\otimes b)]
 -[X\,\iota_0(z\otimes b)\,\iota_0(y\otimes a]
 &=-\langle u\vert a\rangle \bil(y,z)b+\langle u\vert b\rangle\bil(y,z)a\\
 &=\bil(y,z)u=\bil(y,z)X.
\end{split}
\]
Since $\bigl(X\vert\iota_0(y\otimes a)\bigr)=0=\bigl(X\vert \iota_0(z\otimes b)\bigr)$, the result follows in this case.

\smallskip
\item For $X=\iota_0(x\otimes u$, $x\in S_8$, $u\in V$, we have
\[
\begin{split}
[\iota_0(x\otimes u)\,&\iota_0(y\otimes a)\,\iota_0(z\otimes b)]\\
 &=[\langle u\vert a\rangle t_{x,y}+\bil(x,y)t_{u,a},\iota_0(z\otimes b)]\\
 &=\langle u\vert a\rangle\iota_0(\sigma_{x,y}(z)\otimes b)+\bil(x,y)\iota_0(z\otimes \sigma_{u,a}(b))\\
 &=\langle u\vert a\rangle \iota_0\bigl(b\otimes (\bil(x,z)y-\bil(y,z)x)\bigr)-\bil(x,y)\iota_0\bigl(z\otimes(\langle u\vert b\rangle a+\langle a\vert b\rangle u)\bigr).
\end{split}
\]
Thus,
\[
\begin{split}
[X\,\iota_0(y\otimes a)\,&\iota_0(z\otimes b)]-[X\,\iota_0(z\otimes b)\,\iota_0(y\otimes a)]\\
 &=-\bil(x,y)\iota_0\bigr(z\otimes (\langle u\vert b\rangle a+\langle a\vert b\rangle u+\langle u\vert b\rangle a)\bigr)\\
 &\qquad +\bil(x,z)\iota_0\bigl(y\otimes (\langle u\vert a\rangle b+\langle b\vert a\rangle u+\langle u\vert a\rangle b)\bigr)\\
 &\qquad +\bil(y,z)\iota_0\bigl(x\otimes (-\langle u\vert a\rangle b+\langle u\vert b\rangle a)\bigr)\\
 &=\bil(x,y)\langle u\vert a\rangle\iota_0(z\otimes b)-\bil(x,z)\langle u\vert b\rangle \iota_0(y\otimes a)+\bil(y,z)\iota_0(x\otimes u)\\
 &=\bigl(X\vert \iota_0(y\otimes a)\bigr)\iota_0(z\otimes b)-
  \bigl(X\vert\iota_0(z\otimes b)\bigr)\iota_0(y\otimes a)+\bil(y,z)X.
\end{split}
\]

\smallskip
\item For $X=\iota_1(x\otimes 1)$, we have,
\[
\begin{split}
[\iota_1(x\otimes 1)\,&\iota_0(y\otimes a)\,\iota_0(z\otimes b)]\\
 &=[\iota_2(y\bullet x\otimes a),\iota_0(z\otimes b)]\\
 &=\iota_1((y\bullet x)\bullet z\otimes 1)\qquad\textrm{(as $a\bullet 1=-a$ and $a\bullet b=1$),}\\[4pt]
[\iota_1(x\otimes 1)&\iota_0(z\otimes b)\iota_0(y\otimes a)]\\
 &=[\iota_2(z\bullet x\otimes b),\iota_0(y\otimes a)]\\
 &=-\iota_1((z\bullet x)\bullet y\otimes 1)\qquad\textrm{(as $b\bullet 1=-b$, $b\bullet a= -1$).}
\end{split}
\]
Thus,
\[
\begin{split}
[X\,\iota_0(y\otimes a)\,\iota_0(z\otimes b)]&-[X\,\iota_0(z\otimes b)\,\iota_0(y\otimes a)]\\
 &=\iota_1\bigl(((y\bullet x)\bullet z+(z\bullet x)\bullet y)\otimes 1\bigr)\\
 &=\bil(y,z)X,
\end{split}
\]
because the associativity of the bilinear form $\bil$ in a symmetric composition algebra is equivalent to the condition $(x\bullet y)\bullet x=q(x)y=x\bullet (y\bullet x)$ (see \cite[(34.1)]{KMRT}) and hence it follows that $(y\bullet x)\bullet z+(z\bullet x)\bullet y=\bil(y,z)x$ by linearization.

\smallskip
\item For $X=\iota_2(x\otimes 1)$ the situation is similar.

\end{itemize}

Therefore, $\bigl(T,[...],(.\vert .)\bigr)$ is an orthosymplectic triple system and, by its own construction, its Lie superalgebra of inner derivations is isomorphic to $\frel(5;3)$, as $[TT.]=\ad_{[T,T]}=\ad_{\frel(5;3)}$. Thus, the Lie superalgebra $\frtg(T)$ in Proposition \ref{pr:gtOSTS} is isomorphic to the Lie superalgebra $\frg(S_8,S_{1,2})$.

But any derivation $d\in \der T$ extends to a derivation of $\frtg(T)$ which, by Lemma \ref{le:derg83}, is inner. It follows that $\der T=\inder T$ is isomorphic to $\frel(5;3)$, as required.

\smallskip

The associated Lie superalgebra (see Proposition \ref{pr:gOSTS}) is
\begin{equation}\label{eq:gTel53}
\frg=\bigl(\frsp(V)\oplus\frel(5;3)\bigr)\oplus\bigl(V\otimes T\bigr).
\end{equation}
Consider again the description of $\frg(S_8,S_{1,2})$ in \eqref{eq:g83sigmas}:
\[
\frg(S_8,S_{1,2})=\oplus_{\sigma\in\calS_{8,3}}V(\sigma).
\]
Then, as in \eqref{eq:g830sigmas},
\[
\frel(5,3)=\oplus_{\sigma\in\calS_+}V(\sigma),\qquad T=\oplus_{\sigma\in\calS_-}V(\sigma),
\]
with
\[
\begin{split}
\calS_+&=\bigl\{\emptyset,\{1,2,3,4\},\{1,2\},\{3,4\},
           \{2,3,5\},\{1,4,5\},\{1,3,5\},\{2,4,5\}\bigr\},\\[4pt]
\calS_-&=\bigl\{\{5\},\{1,2,5\},\{3,4,5\},\{2,3\},\{1,4\},\{1,3\},\{2,4\}\bigr\}.
\end{split}
\]
Now, assign the index $6$ to the new copy of $V$ in \eqref{eq:gTel53}. Then,
\[
\frg=\oplus_{\sigma\in\tilde\calS}V(\sigma),
\]
with $\tilde\calS\subseteq 2^{\{1,2,3,4,5,6\}}$ given by:
\[
\begin{split}
\tilde\calS&=\bigl\{\emptyset,\{1,2,3,4\},\{1,2\},\{3,4\},
   \{2,3,5\},\{1,4,5\},\{1,3,5\},\{2,4,5\}\\
   &\qquad \{5,6\},\{1,2,5,6\},\{3,4,5,6\},\{2,3,6\},\{1,4,6\},\{1,3,6\},\{2,4,6\}\bigr\}.
\end{split}
\]
Write now
\[
\bar 1=1,\quad \bar 2=3,\quad, \bar 3=5,\quad \bar 4=2,\quad\bar 5=4,\quad \bar 6=6.
\]
Then, we obtain,
\[
\begin{split}
\tilde\calS&=\bigl\{\emptyset,\{\bar 1,\bar 2,\bar 3,\bar 5\},
 \{\bar 1,\bar 4\},\{\bar 2,\bar 5\},
   \{\bar 2,\bar 3,\bar 4\},\{\bar 1,\bar 3,\bar 5\},
   \{\bar 1,\bar 2,\bar 3\},\{\bar 3,\bar 4,\bar 5\}\\
   &\qquad \{\bar 3,\bar 6\},\{\bar 1,\bar 3,\bar 4,\bar 6\},
   \{\bar 2,\bar 3,\bar 5,\bar 6\},
   \{\bar 2,\bar 4,\bar 6\},\{\bar 1,\bar 5,\bar 6\},
   \{\bar 1,\bar 2,\bar 6\},\{\bar 4,\bar 5,\bar 6\}\bigr\},
\end{split}
\]
and this coincides with $\calS_{S_{4,2},S_{4,2}}$ in \cite[\S 5.4]{CE1}. Hence this superalgebra is a Lie superalgebra with the same Cartan matrix $A_{S_{4,2},S_{4,2}}$ in \cite[\S 5.4]{CE1}, thus proving that $\frg$ is isomorphic to the Lie superalgebra $\frg(S_{4,2},S_{4,2})$ in the Supermagic Square.
\end{proof}

\begin{remark}
The previous Theorem shows that the Lie superalgebra $\frel(5;3)$ lives inside $\frg(S_{4,2},S_{4,2})$, and that, in fact, $\frg(S_{4,2},S_{4,2})$ contains a maximal subalgebra isomorphic to $\frsl_2\oplus\frel(5;3)$.
\end{remark}

\bigskip

\section{The Lie superalgebra $\frbr(2;5)$}

In this section  a model of the simple Lie superalgebra $\frbr(2;5)$ is explicitly built. To this aim, consider the $\bZ_2\times\bZ_2$-graded vector space
\[
\frg=\frg_{(0,0)}\oplus\frg_{(1,0)}\oplus\frg_{(0,1)}\oplus\frg_{(1,1)},
\]
with
\[
\begin{split}
\frg_{(0,0)}&=\frsp(V_1)\oplus\frsp(V_2),\\
\frg_{(1,0)}&=\frsp(V_1)\otimes V_2,\\
\frg_{(0,1)}&=V_1\otimes \frsp(V_2),\\
\frg_{(1,1)}&=V_1\otimes V_2,
\end{split}
\]
where, as usual, $V_1$ and $V_2$ are two-dimensional vector spaces endowed with nonzero alternating bilinear forms denoted by $\langle.\vert.\rangle$.

This vector space becomes a superspace with
\[
\frg\subo=\frg_{(0,0)}\oplus\frg_{(1,1)},\quad \frg\subuno=\frg_{(1,0)}\oplus\frg_{(0,1)}.
\]

Now, define a superanticommutative product on $\frg$ by means of the natural Lie bracket on $\frg_{(0,0)}$, the natural action of $\frg_{(0,0)}$ on each $\frg_{(i,j)}$ ($V_i$ is the natural module for $\frsp(V_i)$, while $\frsp(V_i)$ is its adjoint module), and by:
\begin{equation}\label{eq:br25}
\begin{split}
[f\otimes u,g\otimes v]&=\langle u\vert v\rangle[f,g]+2\tr(fg)\gamma_{u,v},\\
[a\otimes p,b\otimes q]&=-\bigl(2\tr(pq)\gamma_{a,b}+\langle a\vert b\rangle [p,q]\bigr),\\
[a\otimes u,b\otimes v]&=\langle u\vert v\rangle\gamma_{a,b}+\langle a\vert b\rangle\gamma_{u,v},\\
[f\otimes u,a\otimes p]&=f(a)\otimes p(u),\\
[f\otimes u,a\otimes v]&=f(a)\otimes \gamma_{u,v},\\
[a\otimes p,b\otimes v]&=-\gamma_{a,b}\otimes p(v),
\end{split}
\end{equation}
for any $a,b\in V_1$, $u,v\in V_2$, $f,g\in\frsp(V_1)=\frsl(V_1)$ and $p,q\in \frsp(V_2)=\frsl(V_2)$. Here, as before, $\gamma_{a,b}=\langle a\vert .\rangle+\langle b\vert .\rangle a$ and similarly for $\gamma_{u,v}$.

This multiplication converts $\frg$ into a $\bZ_2\times\bZ_2$-graded anticommutative superalgebra.

\begin{theorem}
Let $k$ be a field of characteristic $5$. Then the superalgebra $\frg$ above is a Lie superalgebra isomorphic to $\frbr(2;5)$.
\end{theorem}
\begin{proof}
It is clear that all the products in \eqref{eq:br25} are invariant under the action of $\frsp(V_1)\oplus\frsp(V_2)$. Several instances of the Jacobi identity have to be checked. To do so, it is harmless to assume that the ground field $k$ is infinite (extend scalars otherwise) and hence, Zariski topology can be used.

First, for elements $a,b,c\in V_1$ and $u,v,w\in V_2$, to check that the Jacobian
\[
\begin{split}
J(a\otimes u&,b\otimes v,c\otimes w)\\
 &=[[a\otimes u,b\otimes v],c\otimes w]+[[b\otimes v,c\otimes w],a\otimes u]+
    [[c\otimes w,a\otimes u],b\otimes v]
\end{split}
\]
is $0$, it can be assumed, by Zariski density, that $\langle a\vert b\rangle \ne 0$ and $\langle u\vert v\rangle \ne 0$. (Note that the set $\{(a,b)\in V\times V: \langle a\vert b\rangle\ne 0\}$ is a nonempty open set in the Zariski topology of $V\times V$, and hence it is dense.) Moreover, scaling now $b$ and $v$ if necessary, it can be assumed that $\langle a\vert b\rangle =1=\langle u\vert v\rangle$; that is, $\{a,b\}$ is a symplectic basis of $V_1$ and $\{u,v\}$ is a symplectic basis of $V_2$. Now $c=\alpha a+\beta b$ and $w=\mu u+\nu v$ for some $\alpha,\beta,\mu,\nu\in k$. Then:
\[
\begin{split}
J(a\otimes u&,b\otimes v,c\otimes w)\\
 &=[[a\otimes u,b\otimes v],c\otimes w]+[[b\otimes v,c\otimes w],a\otimes u]+
    [[c\otimes w,a\otimes u],b\otimes v]\\
 &=\langle u\vert v\rangle \gamma_{a,b}(c)\otimes w+\langle a\vert b\rangle c\otimes \gamma_{u,v}(w)\\
 &\qquad\qquad +\langle v\vert w\rangle \gamma_{b,c}(a)\otimes u+
    \langle b\vert c\rangle a\otimes \gamma_{v,w}(u)\\
 &\qquad\qquad +\langle w\vert u\rangle \gamma_{c,a}(b)\otimes v+
    \langle c\vert a\rangle b\otimes \gamma_{w,u}(v)\\
 &=(\beta b-\alpha a)\otimes (\mu u+\nu v)+
      (\alpha a+\beta b)\otimes (\nu v-\mu u)\\
 &\qquad\qquad +\mu(\alpha a+2\beta b)\otimes u +
             \alpha a\otimes (\mu u+2\nu v)\\
 &\qquad\qquad -\nu(2\alpha a+\beta b)\otimes v -
             \beta b\otimes(2\mu u+\nu v)\\
 &=0.
\end{split}
\]
Hence, $\frg\subo=\frg_{(0,0)}\oplus \frg_{(1,1)}$ is a Lie algebra, which can be easily checked to be isomorphic to the symplectic Lie algebra $\frsp(V_1\perp V_2)\simeq \frsp_4$.

Now, for elements $f,g,h\in \frsp(V_1)$ and $u,v,w\in V_2$, it can be assumed as before that $\langle u\vert v\rangle =1$ and that $w=\mu u+\nu v$. Then,
\[
\begin{split}
[[f\otimes u,g\otimes v],h\otimes w]
 &=[\langle u\vert v\rangle[f,g]+2\tr(fg)\gamma_{u,v},h\otimes w]\\
 &=\langle u\vert v\rangle [[f,g],h]\otimes w +2\tr(fg)h\otimes \gamma_{u,v}(w),
\end{split}
\]
so that
\begin{equation}\label{eq:jacobian1}
\begin{split}
J(f\otimes u&,g\otimes v,h\otimes w)\\
 &=[[f\otimes u,g\otimes v],h\otimes w]+
    [[g\otimes v,h\otimes w],f\otimes u]+
    [[h\otimes w,f\otimes u],g\otimes w]\\
 &=\langle u\vert v\rangle [[f,g],h]\otimes w
           +2\tr(fg)h\otimes \gamma_{u,v}(w)\\
 &\qquad\qquad +\langle v\vert w\rangle [[g,h],f]\otimes u
           +2\tr(gh)f\otimes \gamma_{v,w}(u)\\
 &\qquad\qquad +\langle w\vert u\rangle [[h,f],g]\otimes v
           +2\tr(hf)g\otimes \gamma_{w,u}(v)\\
 &=\mu\bigl([[f,g],h]-2\tr(fg)h-[[g,h],f]-2\tr(gh)f+4\tr(hf)g\bigr)\otimes u\\
 &\qquad\qquad +\nu\bigl([[f,g],h]+2\tr(fg)h-2\tr(gh)f-[[h,f],g]+2\tr(hf)g\bigr)\otimes v.
\end{split}
\end{equation}
But for any $f\in \frsl(V_1)$, $f^2=-\det(f)=\frac{1}{2}\tr(f^2)id$. Hence $fg+gf=\tr(fg)id$ for any $f,g\in \frsl(V_1)$, and thus
\[
fgf=\tr(fg)f-gf^2=\tr(fg)f-\frac{1}{2}\tr(f^2)g.
\]
Hence,
\[[[g,f],f]=gf^2+f^2g-2fgf=2\tr(f^2)g-2\tr(fg)f
\]
and
\[
[[g,f],h]+[[g,h],f]=4\tr(fh)g-2\tr(fg)h-2\tr(gh)f.
\]
This shows that the Jacobian in \eqref{eq:jacobian1} is trivial. Therefore, the subspace $\frg_{(0,0)}\oplus\frg_{(1,0)}$ is a Lie superalgebra. The same happens to the subspace $\frg_{(0,0)}\oplus \frg_{(0,1)}$.

Take now elements $a,b\in V_1$, $f\in \frsp(V_1)$ and $u,v,w\in V_2$. Then it can be assumed that $\langle a\vert b\rangle =1=\langle u\vert v\rangle$, and $w=\mu u+\nu v$. In this situation:
\[
\begin{split}
J(a\otimes u&,b\otimes v,f\otimes w)\\
 &=[[a\otimes u,b\otimes v],f\otimes w]+
    [[b\otimes v,f\otimes w],a\otimes u]+
   [[f\otimes w,a\otimes u],b\otimes v]
   \\
 &=\langle u\vert v\rangle[\gamma_{a,b},f]\otimes w+\langle a\vert b\rangle f\otimes \gamma_{u,v}(w)\\
 &\qquad\qquad +\gamma_{a,f(b)}\otimes \gamma_{v,w}(u)
     -\gamma_{f(a),b}\otimes \gamma_{u,w}(v)\\
 &=\mu\bigl([\gamma_{a,b},f]-f-\gamma_{a,f(b)}
          -2\gamma_{f(a),b}\bigr)\otimes u\\
 &\qquad\qquad +\nu\bigl([\gamma_{a,b},f]+f-2\gamma_{a,f(b)}
               -\gamma_{f(a),b}\bigr)\otimes v.
 \end{split}
\]
But, since the bilinear map $(c,d)\mapsto \gamma_{c,d}$ is $\frsp(V_1)$-invariant, $[f,\gamma_{a,b}]=\gamma_{f(a),b}+\gamma_{a,f(b)}$.
Hence,
\begin{equation}\label{eq:jacobian2}
J(a\otimes u,b\otimes v,f\otimes w)
 =-\mu\bigl(f+3\gamma_{f(a),b}+2\gamma_{a,f(b)}\bigr)\otimes u
  +\nu\bigl(f-2\gamma_{f(a),b}-3\gamma_{a,f(b)}\bigr).
\end{equation}
Also, by taking the coordinate matrix of $f$ in the symplectic basis $\{a,b\}$, it is checked at once that $f=-\frac{1}{2}\gamma_{f(a),b}+\frac{1}{2}\gamma_{a,f(b)}$. Since the characteristic of $k$ is equal to $5$, this proves that the Jacobian in \eqref{eq:jacobian2} is trivial.

The other instances of the Jacobi identity are checked in a similar way.

\smallskip

Finally, fix symplectic bases $\{a_i,b_i\}$ of $V_i$ ($i=1,2$). Then $\frg$ is $\bZ\times \bZ$-graded by assigning $\deg(a_i)=\epsilon_i$, $\deg(b_i)=-\epsilon_i$, where $\{\epsilon_1,\epsilon_2\}$ denotes the canonical $\bZ$-basis of $\bZ\times\bZ$. Let $\{h_i,e_i,f_i\}$ be the basis of $\frsp(V_i)$ defined as in \eqref{eq:hef}. Then $\espan{h_1,h_2}$ is a Cartan subalgebra of $\frg$, and coincides with the $(0,0)$-component in the $\bZ\times\bZ$-grading. The set of nonzero degrees is
\[
\Phi=\{\pm 2\epsilon_1,\pm 2\epsilon_2,\pm\epsilon_1\pm\epsilon_2,\pm\epsilon_2,\pm 2\epsilon_1\pm\epsilon_2,\pm\epsilon_1, \pm\epsilon_1\pm 2\epsilon_2\}.
\]
Consider the elements
\begin{alignat*}{3}
E_1&=a_1\otimes f_2,&\qquad F_1&=-b_1\otimes e_2,&\qquad H_1&=[E_1,F_1]=-2h_1-h_2,\\
E_2&=h_1\otimes a_2,& F_2&=h_1\otimes b_2,& H_2&=[E_2,F_2]=h_2.
\end{alignat*}
Then, $\espan{H_1,H_2}$ coincides with the previous Cartan subalgebra $\espan{h_1,h_2}$ of $\frg$, $E_1$ belongs to the homogeneous component $\frg_{\epsilon_1-2\epsilon_2}$ in the $\bZ\times\bZ$-grading, and similarly $F_1\in\frg_{-\epsilon_1+2\epsilon_2}$, $E_2\in \frg_{\epsilon_2}$, and $F_2\in\frg_{-\epsilon_2}$. The elements $E_1,E_2,F_1,F_2$ generate the Lie superalgebra $\frg$. Besides,
\[
\begin{split}
[H_1,E_1]&=-2h_1(a_1)\otimes f_2-a_1\otimes [h_2,f_2]
      =-2a_1\otimes f_2+2a_1\otimes f_2=0,\\
[H_1,E_2]&=h_1\otimes (-h_2)(a_2)=-h_1\otimes a_2,\\
[H_2,E_1]&=a_1\otimes [h_2,f_2]=-2a_1\otimes f_2,\\
[H_2,E_2]&=h_1\otimes h_2(a_2)=h_1\otimes a_2,
\end{split}
\]
and similarly for the action of the $H_i$'s on the $F_j$'s. It follows, with the same arguments as in \cite[\S 4]{CE1}, that $\frg$ is the Lie superalgebra with Cartan matrix $\left(\begin{smallmatrix} 0&-1\\ -2&1\end{smallmatrix}\right)$, which is the first Cartan matrix of the Lie superalgebra $\frbr(2;5)$ given in \cite[\S 12]{BGLCartan}.
\end{proof}

\bigskip


\begin{thebibliography}{KMRT98}

\bibitem[BGLa]{BGL}
Sofiane Bouarroudj, Pavel Grozman and Dimitry Leites, \emph{Cartan
matrices and presentations of Cunha and Elduque Superalgebras},
arXiv.math.RT/0611391.

\bibitem[BGLb]{BGLCartan}
\bysame, \emph{Classification of simple finite dimensional modular Lie superalgebras with Cartan matrix},
arXiv:0710.5149 [math.RT].

\bibitem[Bour02]{Bourbaki}
Nicolas Bourbaki, \emph{Lie groups and {L}ie algebras. {C}hapters 4--6}
(Translated from the 1968 French original by Andrew Pressley), Springer-Verlag, Berlin, 2002.

\bibitem[CE07a]{CE1}
Isabel Cunha and Alberto Elduque, \emph{{An extended Freudenthal
Magic Square in characteristic 3}}, J.~Algebra \textbf{317} (2007), 471--509.

\bibitem[CE07b]{CE2}
\bysame, \emph{The extended Freudenthal magic square and Jordan algebras},  Manuscripta Math.  \textbf{123}  (2007),  no.~3, 325--351.

\bibitem[Eld04]{ElduqueMagicI}
Alberto Elduque, \emph{{The magic square and symmetric compositions.}},
Rev. Mat.
  Iberoam. \textbf{20} (2004), no.~2, 475--491.

\bibitem[Eld06a]{NewLook}
\bysame, \emph{A new look at {F}reudenthal's magic square},
Non-associative algebra and its applications, L. Sabinin,
L.V.~Sbitneva, and I.~P. Shestakov, eds., Lect. Notes Pure Appl.
Math., vol. 246, Chapman \& Hall/CRC, Boca Raton, FL, 2006, pp.~149--165. 

\bibitem[Eld06b]{EldNew}
\bysame, \emph{New simple {L}ie superalgebras in characteristic 3},
J. Algebra
  \textbf{296} (2006), no.~1, 196--233. 

\bibitem[Eld07a]{ElduqueMagicII}
\bysame, \emph{{The Magic Square and Symmetric Compositions II}},
Rev. Mat. Iberoamericana \textbf{23} (2007), 57--84.


\bibitem[Eld07b]{EldModular}
\bysame, \emph{{Some new simple modular Lie superalgebras}}, Pacific
J. Math. \textbf{231} (2007), no.~2, 337--359.

\bibitem[Eld]{EldTits3}
\bysame, \emph{The Tits construction and some simple Lie superalgebras in  characteristic 3},  arXiv:math.RA/0703784.

\bibitem[EO02]{EldOkuCompoSuper}
Alberto Elduque and Susumu Okubo, \emph{Composition superalgebras},
Comm.
  Algebra \textbf{30} (2002), no.~11, 5447--5471. 

\bibitem[Jac68]{JacobsonJordan}
Nathan Jacobson, \emph{Structure and representations of {J}ordan
algebras},
  American Mathematical Society Colloquium Publications, Vol. XXXIX, American
  Mathematical Society, Providence, R.I., 1968. 


\bibitem[Kac77]{Kac-Lie}
Victor G. Kac, \emph{Lie superalgebras}, Advances in Math. \textbf{26}
(1977),
  no.~1, 8--96. 

\bibitem[KMRT98]{KMRT}
Max-Albert Knus, Alexander Merkurjev, Markus Rost, and Jean-Pierre
Tignol,
  \emph{The book of involutions}, American Mathematical Society Colloquium
  Publications, vol.~44, American Mathematical Society, Providence, RI, 1998.

\bibitem[Kos70]{Kostrikin}
Alexei I. Kostrikin, \emph{A parametric family of simple {L}ie
algebras}, Izv.
  Akad. Nauk SSSR Ser. Mat. \textbf{34} (1970), 744--756.

\bibitem[McC04]{McCrimmon}
Kevin McCrimmon, \emph{A taste of {J}ordan algebras},
Universitext,
  Springer-Verlag, New York, 2004. 


\bibitem[Oku93]{OkuOTS}
Susumu Okubo, \emph{Triple products and {Y}ang-{B}axter equation.
{I}.
  {O}ctonionic and quaternionic triple systems}, J. Math. Phys. \textbf{34}
  (1993), no.~7, 3273--3291. 


\bibitem[She97]{She97}
Ivan P. Shestakov, \emph{Prime alternative superalgebras of arbitrary
  characteristic}, Algebra i Logika \textbf{36} (1997), no.~6, 675--716, 722.

\bibitem[Tit66]{Tits66}
Jacques Tits, \emph{{Alg\`{e}bres alternatives, alg\`{e}bres de Jordan et
  alg\`{e}bres de Lie exceptionnelles. I: Construction}}, Nederl. Akad. Wet.,
  Proc., Ser. A \textbf{69} (1966), 223--237.



\end{thebibliography}

\providecommand{\bysame}{\leavevmode\hbox to3em{\hrulefill}\thinspace}

\end{document}